\documentclass[12pt]{elsarticle}
\usepackage{amsthm, amsmath,amssymb,amsfonts}
\usepackage{graphicx}
\usepackage{color}

\usepackage{psfig}
\usepackage{graphics}
\usepackage{amsmath}
\usepackage{amssymb}
\usepackage{epsfig}
\usepackage{epstopdf}
\usepackage{color}

\newtheorem{theorem}{Theorem}[section]
\newtheorem{proposition}[theorem]{Proposition}
\newtheorem{corollary}[theorem]{Corollary}
\newtheorem{lemma}[theorem]{Lemma}
\newtheorem{remark}[theorem]{Remark}
\newtheorem{definition}[theorem]{Definition}
\newtheorem{example}[theorem]{Example}

\def\para{\vspace{2mm}}
\def\cP{{\mathcal P}}

\def\gcd{{\rm gcd}}

\def\lc{{\rm lc}}
\def\Content{{\rm Content}}

\def\mult{{\rm mult}}

\def\deg{{\rm deg}}
\def\degree{{\rm deg}}

\def\Res{{\rm Res}}

\begin{document}

\begin{frontmatter}

\title{The Limit Point and the T--function}
\author{Angel Blasco and Sonia P\'erez--D\'{\i}az\\
Departamento de F\'{\i}sica y Matem\'aticas \\
        Universidad de Alcal\'a \\
      28871-Alcal\'a de Henares, Madrid, Spain  \\
angel.blasco@uah.es, sonia.perez@uah.es
}

\begin{abstract} Let  $\cP(t)\in {\Bbb K}(t)^{n}$ be a rational parametrization of an algebraic space curve $\cal C$. In this paper, we introduce the notion of  {\it limit point}, $P_L$, of the given parametrization $\mathcal{P}(t)$, and some remarkable properties  of $P_L$ are obtained. In addition, we generalize the results in \cite{MyB-2017} concerning the {\it T--function}, $T(s)$, which is defined by means of a univariate resultant. More precisely, independently on whether the limit point is regular or not, we show that   $T(s)=\prod_{i=1}^n H_{P_i}(s)^{m_i-1}$, where the polynomials $H_{P_i}(s),\,i=1,\ldots,n$ are the {\it  fibre functions} of the singularities  $P_i\in {\cal C}$ of multiplicity $m_i,\,i=1,\ldots,n$. The roots of $H_{P_i}(s)$ are the fibre of $P_i$ for $i=1,\ldots,n$.   Thus, a complete classification of the singularities of a given space curve, via the factorization of a resultant, is obtained.
\end{abstract}


\end{frontmatter}

\section{Introduction}

A given algebraic curve can be represented in different ways, such as implicitly by
defining polynomials, parametrically by rational functions, or locally parametrically by
power series expansions around a point. These representations all have their individual advantages: an implicit representation allows us to easily decide whether a given point
 lies on a given curve, a parametric representation allows us to generate points
of a given curve over the desired coordinate fields, and using power series
expansions one can, for instance, overcome the numerical problems of tracing a curve
through a singularity.

\para

In the last years, important advances in the study and knowledge of a given algebraic variety (which is in general a curve or a surface)  from its parametric representation has been obtained (see \cite{HSW},  \cite{HL97}, \cite{SWP}, etc.). In particular, an essential problem in computer aided geometric design (CAGD)  is the detection of singularities (see e.g. \cite{MyB-2017},  \cite{Chen2008},  \cite{JSC-Perez}, \cite{MACOM} or \cite{Rubio}).  Understanding the singularities of algebraic curves and surfaces is important for understanding their geometry. In fact, a difficult problem in CAGD is the handling of self-intersections, and the theory of singularities of algebraic varieties is potentially a tool for handling this problem.   For instance, once the singularities are located, one can,  use numerical methods to follow curve branches (see e.g. \cite{Bert}).

\para

In \cite{JSC-Perez}, several results in this sense are provided. In particular,  some formulae
for the computation of the multiplicity of a point are presented. These  formulae simply involve the computation of the degree of a
rational function directly determined from the parametrization.   A further analysis on this topic can be found in \cite{MyB-2017} where, using the  direct relation existing  between the cardinality of the fibre of a given point and its multiplicity, it is shown how easily identify the singularities of the curve as those points whose fibre has more than one element. For this purpose, it is introduced the {\it T--function}, $T(s)$, that is a  polynomial which is obtained  from the computation of a univariate resultant and whose factorization provides the {\it  fibre functions} of the different ordinary singularities as well as their corresponding multiplicities. More precisely, it holds that  $T(s)=\prod_{i=1}^n H_{P_i}(s)^{m_i-1}$, where the polynomials $H_{P_i}(s),\,i=1,\ldots,n$ are the {\it  fibre functions} whose roots are the fibre of the ordinary singularities  $P_i\in {\cal C}$ of multiplicity $m_i,\,i=1,\ldots,n$ (see Theorem 3 in \cite{MyB-2017}). Thus, from this result, one may  easily detect and classify the singularities of the curve.   This topic is also addressed by some other authors; see e.g. \cite{Abhy} and  \cite{Buse2012}.

\para

In this paper, we observe that the relation between fibre and multiplicity fails for one point of the curve, the called {\it limit point}. Given a rational space curve $\mathcal{C}$  over an algebraically closed field of characteristic zero $\mathbb{K}$, and a projective parametrization $\mathcal{P}(t)\in\mathbb{P}^n(\mathbb{K}(t))$ of degree $d$, the limit point is defined as $P_L:=\lim_{t\rightarrow\infty}\mathcal{P}(t)/t^d$. The point $P_L$ is on the curve, since $\mathcal{P}(t)\in\mathcal{C}$ for every $t\in\mathbb{K}$ and $\mathcal{C}$ is a closed set. However, it is not ``well--represented'' by the parametrization and  in fact, its fibre is usually empty, i.e., there is no   $t_0\in\mathbb{K}$ such that $\cP(t_0)=P_L$. We say in this case that the limit point is unreachable via the parametrization $\mathcal{P}(t)$. This circumstance involves some difficulties, since the connection between fibre and multiplicity is lost and then, many results based on that connection do not hold.

\para

Every rational parametrization has a limit point.  If $P_L$ is not an affine point or it is a reachable affine point,   $\cP(t)$ is a {\it normal parametrization}. Otherwise,  if $P_L$ is an affine point and it is not reachable via the parametrization,   $\mathcal{P}(t)$ is not normal and $P_L$ is  the {\it critical point} (see Subsection 6.3 in \cite{SWP}). Thus, under a normal parametrization, every affine point of the curve is reachable via $\mathcal{P}(t)$. However, the problem persists, since $P_L$ remains a point of the curve (affine or at infinity) which is not ``well--represented'' by   $\mathcal{P}(t)$ and, thus, its multiplicity is not the cardinality of its fibre. In fact, some important results presented in  \cite{MyB-2017} and \cite{JSC-Perez} hold for every point of the curve but for $P_L$. In particular, the main theorem in \cite{MyB-2017} (Theorem 3) holds only if $P_L$ is not a singularity.

\para

The main goal in this paper is to explore the nature of the limit point and  analyze the relation between its fibre and its multiplicity. As a remarkable result, we show that $P_L$ is reachable via $\mathcal{P}(t)$ only if it is a singular point of the given curve. In addition, we generalize Theorem 3 in \cite{MyB-2017} independently on whether the limit point is regular or not. In this way, we get a result that allows us to easily compute the ordinary singularities of any rational space curve in any dimension.  A natural but more difficult problem is to consider the case of a given algebraic surface defined  by a rational parametrization. In this case, similar results are expected to be provided. We will deal with this problem in a future work.



\para

The structure of the paper is as follows. In Section \ref{S-notacion}, we provide
the notation and some  previous results. In Section \ref{S-pto-limite} the notion of  {\it limit point}, $P_L$, of the given parametrization $\mathcal{P}(t)$ is introduced and some important properties concerning the multiplicity of $P_L$ are obtained. We show these properties with some illustrative examples. In Section \ref{S-tfunct}, we first summarize some previous properties concerning the T--function introduced in  \cite{MyB-2017}. In particular, Theorem \ref{T-tfunct} and Corollary \ref{C-tfunct} that hold under the assumption that the limit point, $P_L$, is regular. The goal of Section \ref{S-tfunct} is to remove this condition and  generalize both results to the case that $P_L$ is a singularity. The proof of this result as well as a previous technical lemma  appear in Section \ref{S-proof}. Finally, in Section  \ref{S-tfunct}, we  also generalize the theorems obtained  to the case of a given parametric space curve in any dimension (see Subsection \ref{Sub-SpaceCurves}). These results are all illustrated with suitable examples.

\section{Notation and previous results}\label{S-notacion}

Let $\mathcal{C}$ be a rational plane curve  over an algebraically closed field of characteristic zero, $\Bbb K$, defined by the projective parametrization
$$\mathcal{P}(t)=(p_1(t):p_2(t):p(t))\in\mathbb{P}^2(\mathbb{K}(t)),$$
where $\gcd(p_1,p_2,p)=1$. We assume that $\mathcal{C}$ is not a line. Let $d_1=\deg(p_1)$,
$d_2=\deg(p_2)$, $d_3=\deg(p)$ and $d=\max\{d_1,d_2,d_3\}$. Then, we   write
$$\left\{\begin{array}{r}p_1(t)=a_0+a_1t+a_2t^2+\cdots+a_{d}t^{d}\\p_2(t)=b_0+b_1t+b_2t^2+\cdots+b_{d}t^{d}\\p(t)=c_0+c_1t+c_2t^2+\cdots+c_{d}t^{d}\end{array}\right.$$

Associated with  ${\mathcal P}(t)$, we
consider the induced rational map $\psi_{\mathcal
P}:{\mathbb{K}}\longrightarrow {\mathcal C} \subset \mathbb{P}^2(\mathbb{K}); t\longmapsto
{\mathcal P}(t)$, and $\degree(\psi_{\mathcal P})$ denotes the degree of the rational map
$\psi_{\mathcal P}$ (see e.g. \cite{Harris:algebraic} pp. 80 or \cite{shafa}
pp. 143). As an important result,
we recall that
 the birationality of  $\psi_\cP$, i.e. the properness of $\cP(t)$, is
characterized  by $\deg(\psi_\cP)=1$ (see \cite{Harris:algebraic}
and \cite{shafa}).  Also, we recall that the degree of a rational
map can be seen as the cardinality of the fibre of a generic element
(see Theorem 7, pp. 76 in \cite{shafa}). We will use this
characterization in our reasoning. For this purpose, we denote by
${\mathcal F}_{\mathcal P}(P)$ the fibre of a point $P\in \mathcal
C$ via the parametrization $\mathcal{P}(t)$; that is $ {\mathcal F}_{\mathcal P}(P)={\mathcal P}^{-1}(P)=\{
t\in \mathbb{K} \,|\, {\mathcal P}(t)=P \}. $

\para

It is well known  that almost all points of $\cal C$ (except at most a finite number of points)
are generated via $\cP(t)$ by the same number of parameter values, and this number is
the degree of $\psi_\cP$ (see Subsection 2.2. in \cite{SWP}). Thus, intuitively speaking, the degree of $\psi_\cP$ measures the number of times that
${\mathcal P}(t)$ traces the curve when the parameter takes values in $\Bbb K$. Taking into
account this intuitive notion, the degree of
the mapping $\psi_\cP$ is also called   the {\it tracing index} of $\cP(t)$. In Chapter 4 in \cite{SWP}, it is proved that the tracing index of   $\mathcal{P}(t)$ can be computed as $\deg(\psi_\mathcal{P}(t))=\deg_t(G)$, where
\begin{equation}\label{Eq-fibra-generica}\left\{\begin{array}{l}G_1(s,t)=p_1(s)p(t)-p(s)p_1(t)\\G_2(s,t)=p_2(s)p(t)-p(s)p_2(t)\\G_3(s,t)=p_1(s)p_2(t)-p_2(s)p_1(t)\end{array}\right.\end{equation}
and $G(s,t)=\gcd(G_1(s,t),G_2(s,t),G_3(s,t))$.

\para

The cardinality of the fibre of    $\psi_\cP$ is the same for almost all
points on $\cal C$; that is,   all
but finitely many points in $\cal C$ are generated, via $\cP(t)$, by exactly $\deg_t(G)$ parameter values.  Nevertheless, for finitely many exceptions, the cardinality may vary. We can compute the fibre of a particular point $P=(a,b,c)$ by solving the corresponding {\it fibre equations}
\begin{equation}\label{Eq-fibre-equations}\left\{\begin{array}{l} \phi_1(t):=ap(t)-cp_1(t)=0\\ \phi_2(t):=bp(t)-cp_2(t)=0\\\phi_3(t):=ap_2(t)-bp_1(t)=0\end{array}\right.\end{equation}
Observe that $\mathcal{P}(t_0)=P$ if and only if $\phi_1(t_0)=\phi_2(t_0)=\phi_3(t_0)=0$. This motivates the following definition.

\para

\begin{definition}\label{D-fibre-func}
Let $P\in\mathbb{P}^2(\mathbb{K})$, and a parametrization
$\mathcal{P}(t)\in\mathbb{P}^2(\mathbb{K}(t))$. We  define the {\it fibre function} of $P$ via
$\mathcal{P}(t)$ as $H_P(t):=\gcd(\phi_1,\phi_2,\phi_3).$
\end{definition}

\begin{remark}\label{R-afin-noafin} Note that the roots of $H_P$ determine the fibre of $P$. 
In addition, we observe  that the above expression for $H_P$ may be simplified if we consider  the following cases (see Remark 2 in \cite{MyB-2017}):
\begin{itemize}
\item If $P$ is an affine point then $H_P(t)=\gcd(\phi_1(t),\phi_2(t)).$
\item If $P$ is a point at infinity then $H_P(t)=\gcd(p(t),\phi_3(t)).$
\end{itemize}
\end{remark}

\para

Throughout this paper we  assume that $\mathcal{P}(t)$ is a proper parametrization; otherwise, we  obtain a proper one by reparametrizing $\cal P$ (see e.g. \cite{Perez-Repara1}). This means that $\deg(\psi_\mathcal{P}(t))=1$ and, so, the cardinality of the fibre is 1 for almost every point of the curve. However, this cardinality may be different for finitely many points. In fact,  the fibre of a singular point is  greater than 1  (see e.g. \cite{JSC-Perez}). The method proposed in \cite{MyB-2017} for computing the ordinary singularities of a rational curve is based on this idea. On the other hand, the cardinality of the fibre may be less than $1$ if we consider the {\it limit point} of the parametrization.

\para

\begin{definition}\label{D-pto-limite}
We define the {\it limit point} of the parametrization $\mathcal{P}(t)$
as
$$P_L=\lim_{t\rightarrow\infty}\mathcal{P}(t)/t^d=(a_d:b_d:c_d).$$
\end{definition}

\para

From this definition, it follows that every parametrization has only one limit point. In addition, $P_L\in {\cal C}$  since $\mathcal{P}(t)/t^d \in\mathcal{C}$ for every $t\in\mathbb{K}$ and $\mathcal{C}$ is a closed set. However, $P_L$  is not a conventional point. The  following results, which are proved in \cite{JSC-Perez}, hold for every point of the curve but for $P_L$.

\para

\begin{theorem}\label{T-tangentes}
Let $\mathcal{C}$ be a rational algebraic curve defined by a proper  pa\-ra\-me\-tri\-za\-tion $\mathcal{P}(t)$, with limit point $P_L$. Let $P\neq P_L$ be a point of $\mathcal{C}$ and let $H_P(t)=\prod_{i=1}^n(t-s_i)^{k_i}$ be its fibre function. Then, $\mathcal{C}$ has  $n$ tangents at $P$ of multiplicities   $k_1,\ldots,k_n$, respectively.
\end{theorem}

\begin{corollary}\label{C-multiplicidad}
Let $\mathcal{C}$ be a rational algebraic curve defined by a proper pa\-ra\-me\-tri\-za\-tion $\mathcal{P}(t)$, with limit point $P_L$. Let $P\neq P_L$ be a point of $\mathcal{C}$ and let $H_P(t)$ be its fibre function. Then, $\mult_P(\mathcal{C})=\deg(H_P(t)).$
\end{corollary}


\para

\begin{figure}[h]
$$
\begin{array}{cc}
\psfig{figure=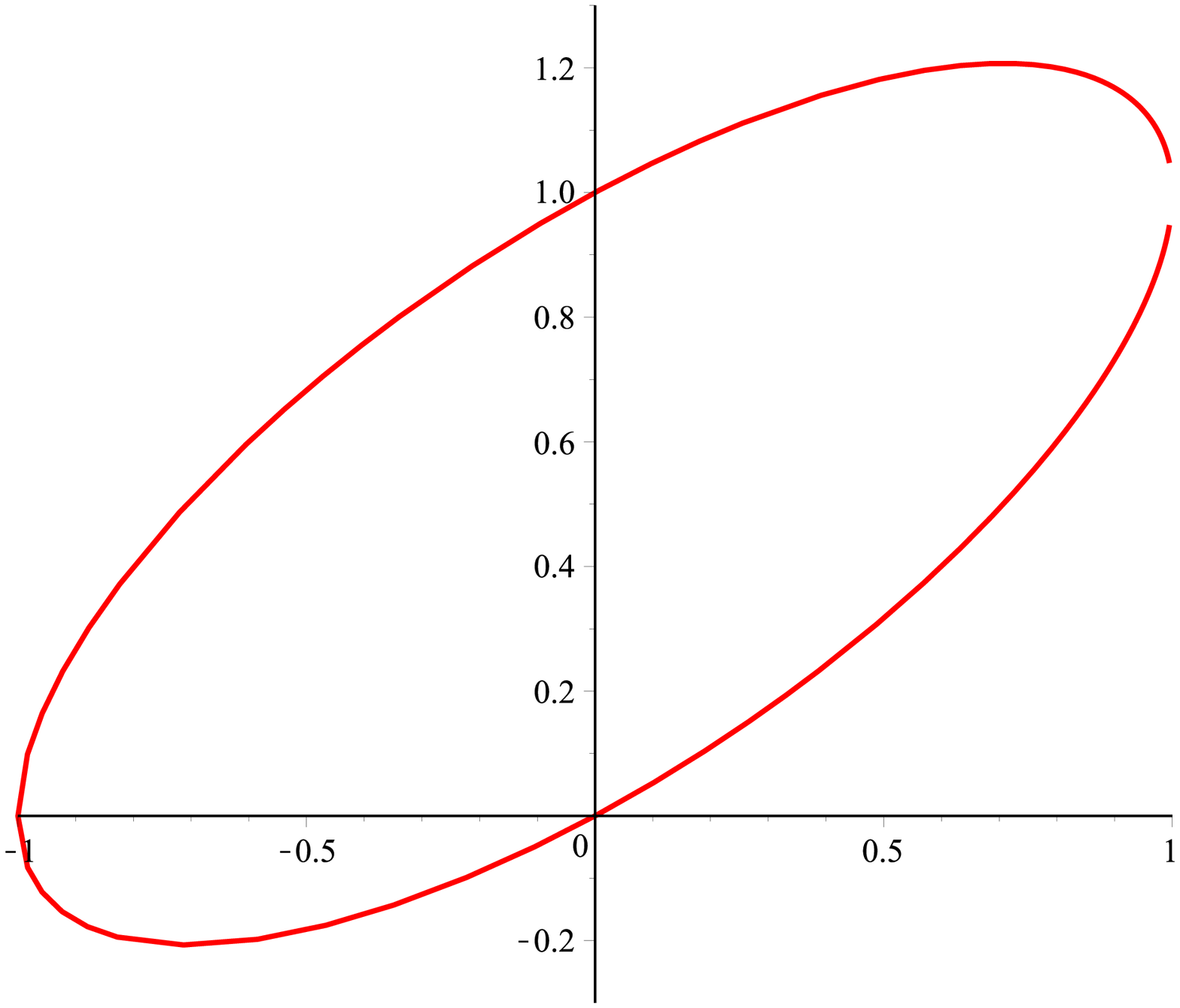,width=4.5cm,height=4.5cm,angle=0} &
\psfig{figure=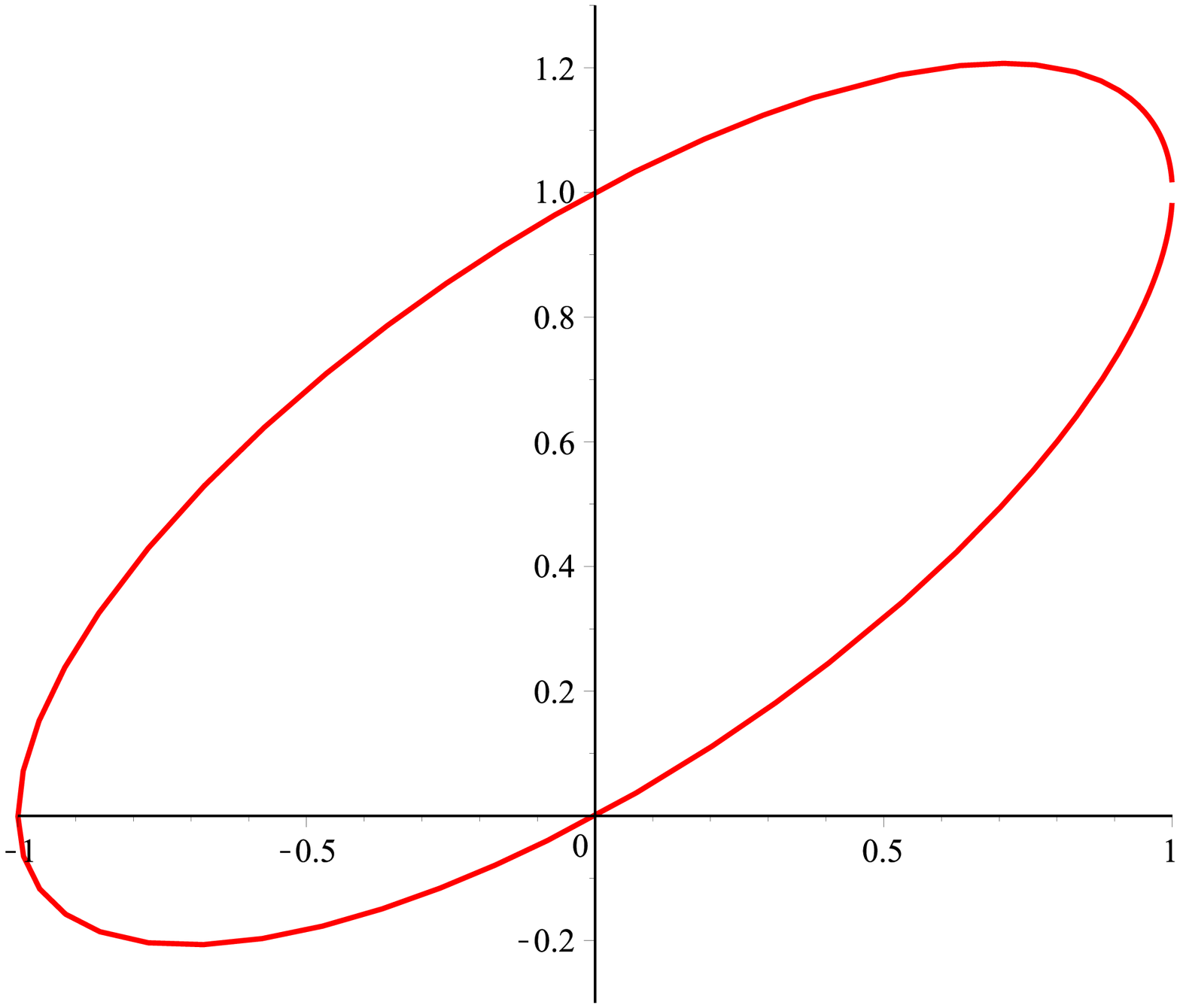,width=4.5cm,height=4.5cm,angle=0}
\end{array}
$$ \caption{Curve $\mathcal{C}$ plotted from $\cP(t)$ with $-20\leq t\leq 20$ (left) and  $-60\leq t\leq 60$ (right)}\label{F-pto-limite-unreach}
\end{figure}

Theorem \ref{T-tangentes} and Corollary \ref{C-multiplicidad} show that there exists a strong relation between the fibre of a point and its multiplicity, but they fail if the point is $P_L$ since this point is not ``well--represented'' by the parametrization. In fact, most of the times it holds that $\mathcal{F}_\mathcal{P}(P_L)=\emptyset$, i.e., there is no $t_0\in\mathbb{K}$  such that $\mathcal{P}(t_0)=P_L$. In this case, we say that the limit point is {\it unreachable} via the parametrization. In order to illustrate the concept of limit point, let us consider the ellipse defined by the projective parametrization
$$\cP(t)=(t^2-1,t^2-t,t^2+1)\in\mathbb{P}^2(\mathbb{C}(t))$$
The limit point is, in this case, the affine point $P_L=(1,1)$. In Figure \ref{F-pto-limite-unreach} (left), we plot the curve using $\cP(t)$ with $-20\leq t\leq 20$. In Figure \ref{F-pto-limite-unreach} (right), we  plot it using $\cP(t)$ with  $-60\leq t\leq 60$. Note that $P_L$ is not reached by the parametrization but it would be reached in the limit, when  $t$ tends to $\infty$.

\para

We say that $\cP(t)$ is a {\it normal parametrization} if $P_L$ is an infinity point or it is a reachable affine point; otherwise, $\cP(t)$ is not normal and it is said that $P_L$ is  the {\it critical point} (see Subsection 6.3 in \cite{SWP}).

\para

In the next section some important properties as well as essential results concerning the limit point are obtained. In particular, we see that $P_L$ is a special point of $\cal C$ whose multiplicity has to be carefully computed. In addition, we prove that if the limit point is reached by $\cP(t)$, then it is a singularity. These results will be used in Section \ref{S-tfunct}, where the relation with the  {\it T--function} is studied. The  T--function provides essential information about the singularities of the given curve $\cal C$. More precisely, its factorization provides the fibre functions of each singularity as well as its corresponding multiplicity (see \cite{MyB-2017}).

\section{The limit point and the hidden multiplicity}\label{S-pto-limite}


In Definition \ref{D-pto-limite}, we introduce the notion of limit point of the parametrization $\mathcal{P}(t)$ as
$$P_L=\lim_{t\rightarrow\infty}\mathcal{P}(t)/t^d=(a_d:b_d:c_d),$$
where $d=\max\{d_1,d_2,d_3\}$, and  $d_1=\deg(p_1)$,
$d_2=\deg(p_2)$, $d_3=\deg(p)$. One may determine the fibre of $P_L$ from the corresponding fibre function, that we denote for this particular point as $H_L(t)$. That is,
$$H_L(t):=\gcd(\phi^L_1,\phi^L_2,\phi^L_3),$$
where
\begin{equation}\label{Eq-fibra-PL}\left\{\begin{array}{l} \phi^L_1(t)=a_dp(t)-c_dp_1(t)\\ \phi^L_2(t)=b_dp(t)-c_dp_2(t)\\\phi^L_3(t)=a_dp_2(t)-b_dp_1(t).\end{array}\right.\end{equation}

\para

The functions $\phi^L_1$, $\phi^L_2$ and $\phi^L_3$ are obtained as particular cases of those introduced in (\ref{Eq-fibre-equations}). Because of the importance of this point, we use a specific notation for them.

\para

\begin{remark}\label{R-afin-noafin-PL}
 Remark \ref{R-afin-noafin} can  be applied similarly for this special case (when the limit point is considered). More precisely,  depending on whether $P_L$ is an affine point or an infinity point, its   fibre function can be expressed as follows:
\begin{itemize}
\item If $P_L$ is an affine point, then $H_L(t)=\gcd(\phi^L_1(t),\phi^L_2(t)).$
\item If $P_L$ is a point at infinity, then $H_L(t)=\gcd(p(t),\phi^L_3(t)).$
\end{itemize}
\end{remark}

\para

Note that we can not compute the multiplicity of $P_L$ from the fibre since Corollary \ref{C-multiplicidad} does not hold in this particular case. Although $P_L$ is a point of $\cal C$, it is not ``well--represented'' by the parametrization $\cP(t)$ and hence, the cardinality of the fibre does not provide its   multiplicity. In the following, and throughout this section, we illustrate this statement and we present a method that allows to compute the multiplicity of $P_L$.

\para

For this purpose, we consider a reparametrization of $\cP(t)$,   $\mathcal{U}(t)$, such that $P_L$ is not the limit point (that is, ${\cal U}(t)$ has a limit point $U_L\not=P_L$). First, we  assume that   $P_L$ is not reachable via $\cP(t)$, i.e., there is no $s_0\in\mathbb{K}$ such that $\mathcal{P}(s_0)=P_L$ (afterwards, we will analyze the case of  $P_L$  reached by the parametrization $\cP(t)$). In this case,  $\mathcal{F}_{\mathcal{P}}(P_L)=\emptyset$ and the system   in  (\ref{Eq-fibra-PL}) does not have any solution. Then, we consider $\mathcal{U}(t)=\mathcal{P}(1/t)$, and we write  $\mathcal{U}(t)=(u_1(t):u_2(t):u(t))$, where
$$\left\{\begin{array}{l}u_1(t)=p_1(1/t)t^d=a_0t^d+a_1 t^{d-1}+\cdots+a_{d}\\u_2(t)=p_2(1/t)t^d=b_0t^d+b_1 t^{d-1}+\cdots+b_{d}\\u(t)=p(1/t)t^d=c_0t^d+c_1 t^{d-1}+\cdots+c_{d}\end{array}\right.$$
The limit point of $\mathcal{U}(t)$ is
$$U_L=\lim_{t\rightarrow\infty}\mathcal{U}(t)/t^d=\left(a_0:b_0:c_0\right)=\mathcal{P}(0)$$
(note that $a_0=b_0=c_0=0$ is not possible since we are assuming that $\gcd(p_1,p_2,p)$=1). On the other hand, $P_L$ is a usual point of $\cal C$, which can be obtained as $P_L=\mathcal{U}(0)$. Since $\mathcal{U}(t)=\mathcal{P}(1/t)$ and $\cP(t)$ is proper, we get that $\mathcal{U}(t)$ is also a proper para\-me\-tri\-zation. Thus, we may apply Corollary \ref{C-multiplicidad} to determine the multiplicity of $P_L$ by computing the cardinality of $\mathcal{F}_{\mathcal{U}}(P_L)$.

\para

In order to get it, we obtain the corresponding fibre function, that we denote as $H^{\cal U}_{L}(t)$. We observe that if $H^{\cal U}_L(s_i)=0$, for some $s_i\neq 0$, then $\mathcal{U}(s_i)=P_L$ and thus $\mathcal{P}(1/s_i)=P_L$, which is impossible by assumption, since  $P_L$ can not be reached by $\cP(t)$. Thus, we have that $H^{\cal U}_{L}(t)=t^r$ for some $r\in {\Bbb N}$. Moreover, we have that $r\geq 1$, since $\mathcal{U}(0)=P_L$.


\para

Therefore, from Corollary  \ref{C-multiplicidad}, we get that the multiplicity of $P_L$ is $r\geq 1$. We refer to this  multiplicity (that can not be obtained from the parametrization  $\mathcal{P}(t)$)  as the {\it hidden multiplicity of the limit point $P_L$}, and we represent it as $m_H$. In the following, we illustrate the above procedure with an example.

\para

\begin{example}\label{Ex-PL-unreachable}
Let $\mathcal{C}$ be a rational plane curve defined by the projective parametrization
$$\mathcal{P}(t)=(t^6+2t^5+2t^4+3t^3+2t^2+t+1:t^6+2t^5+t^4+t^3+t^2:t^4+t^2)\in\mathbb{P}^2(\mathbb{C}(t)).$$
The limit point is  $P_L=\lim_{t\rightarrow\infty}{\cP}(t)/t^6=(1:1:0).$
In order to compute the fibre of $P_L$, we solve the system given in (\ref{Eq-fibra-PL}) that can be expressed as (see Remark \ref{R-afin-noafin-PL})
$$\left\{\begin{array}{l}p(t)=t^4+t^2=0\\p_1(t)-p_2(t)=2t^3+t+1=0.\end{array}\right. $$
This system does not have any solution, which implies that  $P_L$ is not reached by the parametrization $\cP(t)$. We consider the reparametrization
$$\mathcal{U}(t)=\mathcal{P}\left(1/t\right)=(t^6+t^5+2t^4+3t^3+2t^2+2t+1:t^4+t^3+t^2+2t+1:t^4+t^2)$$
and now $P_L$ is a usual point obtained as $P_L=\mathcal{U}(0)$. Furthermore, the fibre function is    $H^{\cal U}_{L}(t)=t^2$ and thus, from Corollary  \ref{C-multiplicidad}, we conclude that $P_L$ is a singular point of $\cal C$ of multiplicity $2$ (note that $m_H=2$).
\end{example}

\para

Now, let us assume that the limit point,  $P_L$, can be reached by the parametrization, i.e. $\mathcal{F}_{\mathcal{P}}(P_L)\neq\emptyset$. Then,
 $H_L(t)=\prod_{i=1}^n(t-s_i)^{k_i}$. If we could apply Corollary  \ref{C-multiplicidad},
we would get that the multiplicity of $P_L$ is
$\deg(H_L(t))=k_1+\cdots+k_n$. However, this is not true since the fibre of $P_L$, via the parametrization $\cP(t)$, does not determine correctly the multiplicity  of $P_L$.

\para

In order to illustrate this  statement, we
reason as above, and we consider the reparametrization $\mathcal{U}(t)=\mathcal{P}(1/t)$. Let us assume that $s_i\not=0$ for $i=1,\ldots,n$ (see Remark \ref{R-raiz-cero}). Then, for each root $s_i$ of $H_L$ we have that $1/s_i$ is a root  of $H^{\cal U}_L$, since $\mathcal{U}(1/s_i)=\mathcal{P}(s_i)=0$. In addition, $H^{\cal U}_L$ has one new root given by $t=0$, since $\mathcal{U}(0)=P_L$. Thus, from $H_L(t)=\prod_{i=1}^n(t-s_i)^{k_i}$, we easily get that
\begin{equation}\label{Eq-HLQ}H^{\cal U}_L(t)=t^r\prod_{i=1}^n(t-1/s_i)^{k_i}.\end{equation}

Finally, by applying Corollary \ref{C-multiplicidad} (note that $P_L$ is not the limit point of ${\cal U}(t)$), we conclude that
$$\mult_{P_L}(\mathcal{C})=\deg(H^{\cal U}_L(t))=r+k_1+\cdots+k_n.$$

Observe that part of this multiplicity is $k_1+\cdots+k_n=\deg(H_L(t))$. It is given by the degree of the fibre function $H_L(t)$, which is obtained from the parametrization $\mathcal{P}(t)$. We refer to this multiplicity  as
the {\it visible multiplicity of the limit point $P_L$} . However, there is another part, $r$, that could not be obtained via $\cP(t)$; this is the hidden multiplicity ($m_H$) introduced above.

\para

\begin{remark}\label{R-raiz-cero} \begin{enumerate}
\item The above reasoning is not correct if $s_i=0$ for some $i=1,\ldots,n$, since  $1/s_i$ is not a root of the polynomial $H^{\cal U}_L(t)$. We observe that in this situation, $P_L$ is the limit point also with the parametrization ${\cal U}(t)$ (note that $U_L=\mathcal{P}(0)=P_L$) and thus, we can not apply Corollary  \ref{C-multiplicidad}. This problem can be solved by considering a new reparametrization of the form $\mathcal{Q}(t)=\mathcal{P}(\theta t/(t-1))$, where $\theta\neq 0$ and $\theta\neq s_i$ for every $i=1,\ldots,n$ (see Section \ref{S-tfunct}). Reasoning with this reparametrization, we get an equivalent result.
\item  Many authors (see e.g. \cite{Chen2008}) use the  homogeneous parametrization  $\overline{\mathcal{P}}(t,h)=(p_1(t/h)h^d:p_2(t/h)h^d:p(t/h)h^d)$. The multiplicity of $P_L$ is  correctly obtained by using this parametrization. Observe that every point of  $\cal C$ is reachable by $\overline{\mathcal{P}}$ since $\mathcal{P}(t)=\overline{\mathcal{P}}(t,1)$ for every $t\in\mathbb{K}$, and  $\overline{\mathcal{P}}(1,h)=\mathcal{U}(h)$ (which implies that $P_L=\overline{\mathcal{P}}(1,0)$). Therefore,  $\overline{\mathcal{P}}(t,h)$ provides the visible and the hidden multiplicity of the limit point $P_L$. In fact, it is easy to check that the corresponding fibre function is
$$\overline{H}_L(t,h)=h^r\prod_{i=1}^n(t-s_ih)^{k_i},$$
where $\cP(s_i)=P_L,\,i=1,\ldots,n,\,s_i\not=s_j,\,i\not=j$, $k_1+\cdots+k_n$ is the visible multiplicity and $r$ is the hidden one (see the analogy with (\ref{Eq-HLQ})).
\end{enumerate}\end{remark}

In the following proposition, we prove an important property concerning the limit point. Namely, if the limit point is reached by $\cP(t)$, then it is a singularity. However, clearly the reciprocal is not true (see Example \ref{Ex-PL-unreachable}).

\para

\begin{proposition}\label{P-alcanzable-singular} If  the limit point is reached by the parametrization $\cP(t)$, then it is a singularity.
\end{proposition}
\noindent\textbf{Proof:} The multiplicity of $P_L$ is given by
$$\mult_{P_L}(\mathcal{C})=k_1+\cdots+k_n+m_H,$$
where $k_1+\cdots+k_n=\deg(H_L(t))$ is the visible multiplicity, and $m_H$ is the hidden multiplicity. Note that $m_H\geq 1$ ($\mathcal{U}(0)=P_L$), and since $P_L$ is reached by the parametrization, then some of its multiplicity is visible, i.e.  $k_1+\cdots+k_n\geq 1$. Therefore, $\mult_{P_L}(\mathcal{C})\geq 2$ and thus, $P_L$ is a singularity.\hfill $\Box$

\para

\begin{example}\label{Ex-PL-reachable} Let  $\mathcal{C}$ be the rational curve defined over $\Bbb C$ by the projective parametrization $\mathcal{P}(t)=$  $$(t^6+2t^5+3t^4+3t^2+4t^3+2t+1: -t^4-t^3-t^2-2t-1: t^7+3t^5+t^4+3t^3+t+2t^2+1).$$
The limit point is
$P_L=\lim_{t\rightarrow\infty}\mathcal{P}(t)/t^6=(0:0:1).$ The fibre of $P_L$ is obtained by solving the system   in (\ref{Eq-fibra-PL}), that can be expressed as (see Remark \ref{R-afin-noafin-PL})
$$\left\{\begin{array}{l}p_1(t)=t^6+2t^5+3t^4+3t^2+4t^3+2t+1=0\\p_2(t)=-t^4-t^3-t^2-2t-1=0.\end{array}\right.$$
The gcd of both polynomials is  $H_{L}(t)=(t+1)$, and thus ${\mathcal F}_{\mathcal P}(P_L)=\{-1\}$ (i.e. $\mathcal{P}(-1)=(0:0:1)=P_L$). If we could apply Corollary \ref{C-multiplicidad}, we would deduce that $P_L$ is a regular point since its multiplicity is $1$. However, this is not true. Indeed: let us consider the reparametrization $\mathcal{U}(t)=\cP(1/t)=$ $$(t^7+2t^6+3t^5+4t^4+3t^3+2t^2+t:-t^7-2t^6-t^5-t^4-t^3:t^7+t^6+2t^5+3t^4+t^3+3t^2+1).$$

\begin{figure}[h]
$$
\begin{array}{cc}
\psfig{figure=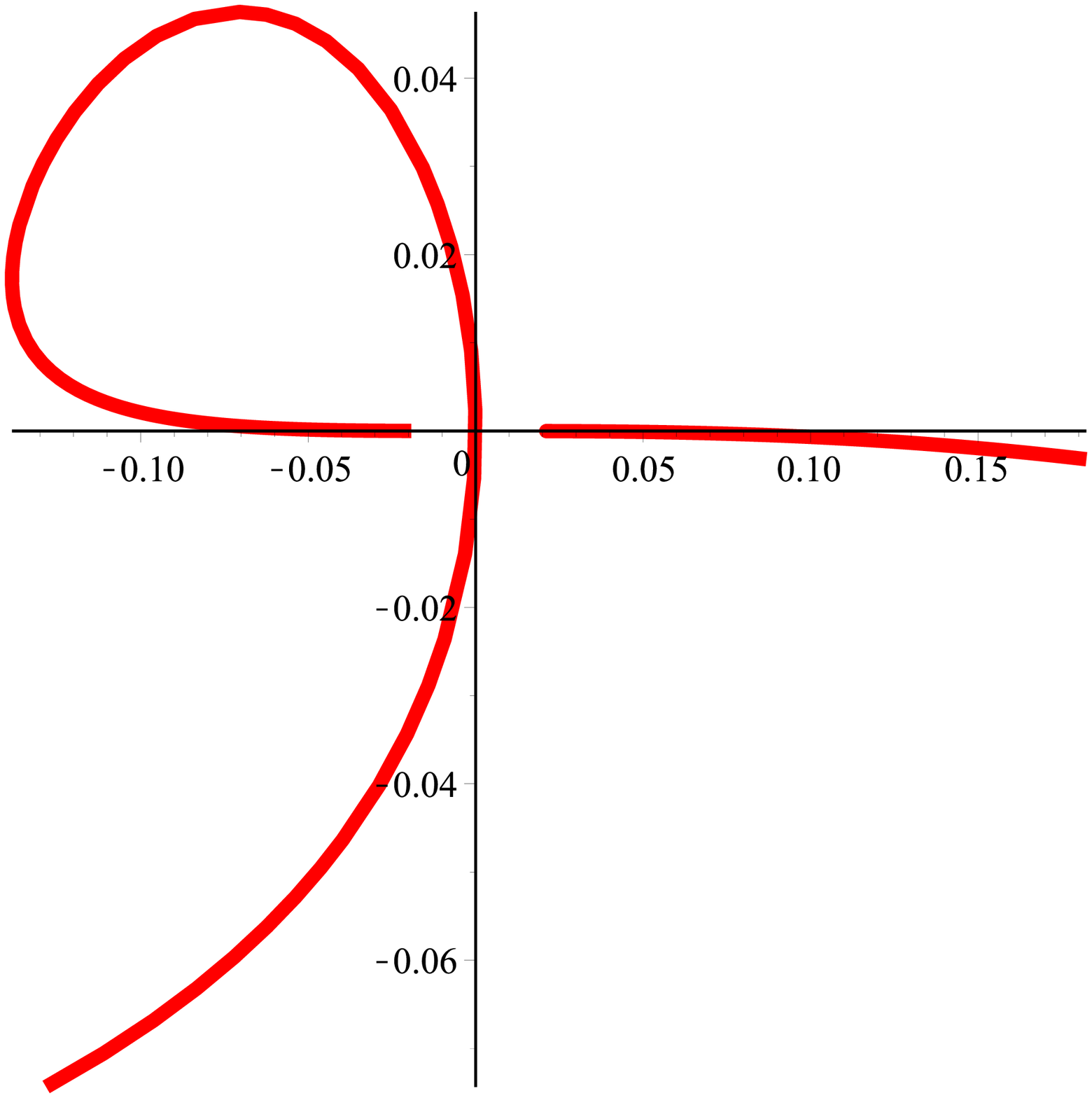,width=4.5cm,height=4.5cm,angle=0} &
\psfig{figure=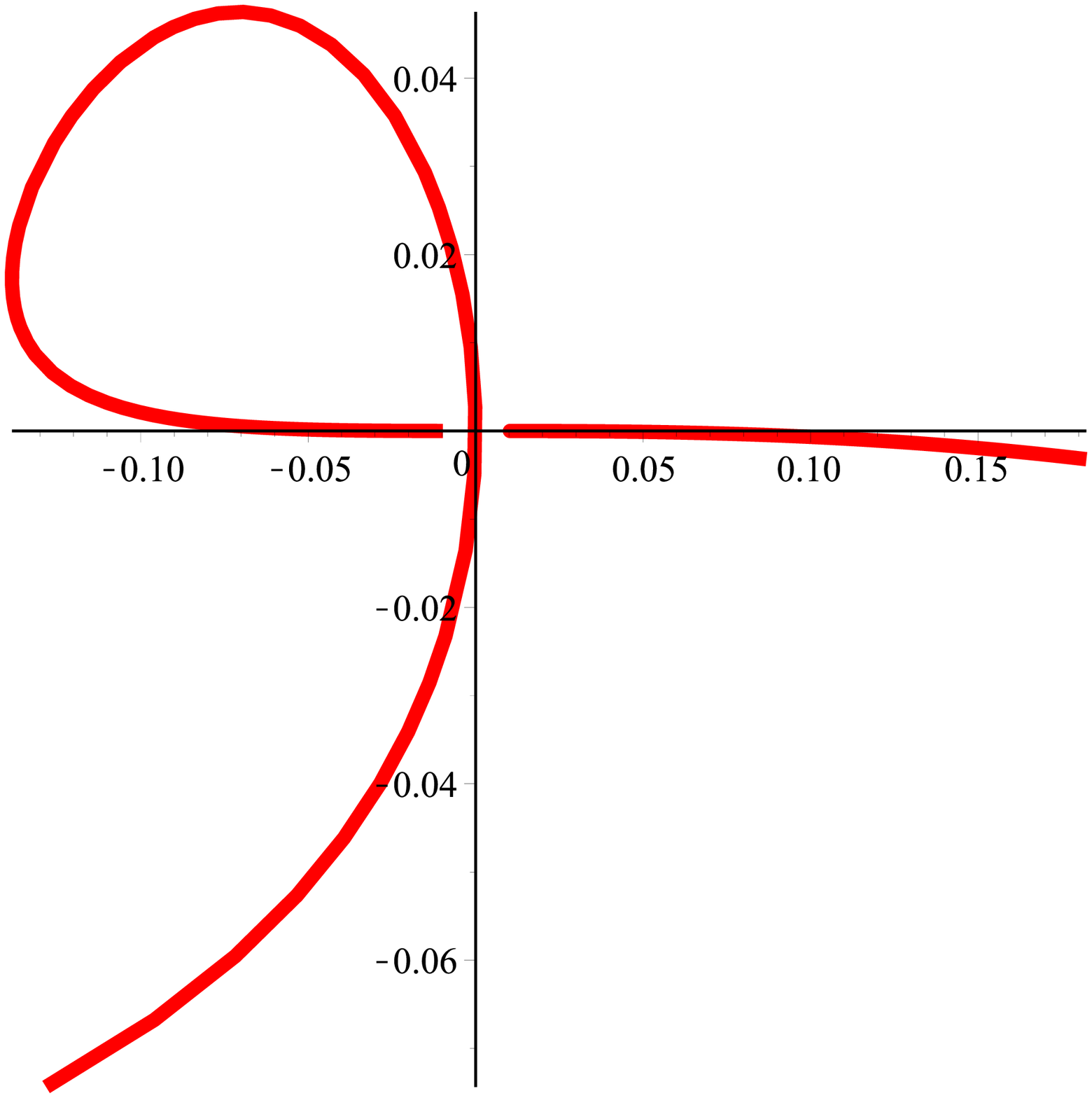,width=4.5cm,height=4.5cm,angle=0}
\end{array}
$$ \caption{Curve $\mathcal{C}$ plotted from $\cP(t)$ with $-50\leq t\leq 50$ (left) and  $-100\leq t\leq 100$ (right)}\label{F-pto-limite0}
\end{figure}

Note that  $\mathcal{U}(0)=(0:0:1)=P_L$. The fibre of $P_L$ (via ${\cal U}(t)$), ${\mathcal F}_{\mathcal U}(P_L)$, is given by the common roots of the equations
$$\left\{\begin{array}{l}q_1(t)=t^7+2t^6+3t^5+4t^4+3t^3+2t^2+t=0\\q_2(t)=-t^7-2t^6-t^5-t^4-t^3=0.\end{array}\right.$$
The gcd of both polynomials is  $H^{\cal U}_{L}(t)=t(t+1)$ and thus, the cardinality of the fibre, ${\mathcal F}_{\mathcal U}(P_L)$, is $2$. Hence, $P_L$ is a double point of $\cal C$. The visible multiplicity is $1$, and the hidden multiplicity is $1$.

\para

 In Figure \ref{F-pto-limite0},  we plot the curve $\mathcal{C}$ in a neighborhood of the limit point $P_L$. More precisely, in Figure \ref{F-pto-limite0} (left), we   plot $\mathcal{C}$ using $\cP(t)$ with $-50\leq t\leq 50$. In Figure \ref{F-pto-limite0} (right), we  plot $\mathcal{C}$ using $\cP(t)$ with  $-100\leq t\leq 100$. Note that  $P_L=(0,0)$ is reached once by $\cP(t)$ although it would be reached again in the limit, when  $t$ tends to $\infty$. This second time (as the limit of the parametrization) does not have any fibre but it provides a second tangent (which increases the multiplicity of $P_L$) that is not detected from the parametrization  $\mathcal{P}(t)$.
\end{example}

\section{The limit point and the T--function}\label{S-tfunct}

In  \cite{MyB-2017}, a method for computing the singularities of a rational algebraic curve from its parametric representation is proposed. The method is based on the construction and factorization of a polynomial called the {\it T--function}, which may be defined in the following three ways:
\begin{equation}\label{Eq-tfunct}T(s)=\frac{R_{12}(s)}{p(s)^{\lambda_{12}-1}}=\frac{R_{13}(s)}{p_1(s)^{\lambda_{13}-1}}=\frac{R_{23}(s)}{p_2(s)^{\lambda_{23}-1}},\end{equation}
where
$$R_{ij}(s)=\Res_t\left(\frac{G_i(s,t)}{t-s},\frac{G_j(s,t)}{t-s}\right)$$
and $\lambda_{ij}=\min\{\deg_t(G_i),\deg_t(G_j)\}$,\,$i,j\in \{1,2,3\},\,i<j$.

\para

The  T--function provides essential information about the singularities of the given curve $\cal C$. More precisely,   the factorization of the T--function gives the fibre functions of the singularities of $\cal C$.  We remark that from the fibre function of a point $P$, one may determine the multiplicity of $P$  as well as its fibre and the tangent lines at $P$. In \cite{MyB-2017}, some important results concerning the T--function are shown. In the following, we summarize some of them.

\para

\begin{lemma}\label{L-tfunct} Let $\mathcal{C}$ be a rational algebraic plane curve defined by a parametrization $\mathcal{P}(t)$, with limit point $P_L$. Let $P\neq P_L$  be an ordinary singular point of multiplicity $m$. It holds that
$$ T(s)=H_P(s)^{m-1}T^*(s),$$
where $T^*(s)\in {\Bbb K}[s]$ and $\gcd(H_P(s),T^*(s))=1.$
\end{lemma}

\begin{theorem}\label{T-tfunct} {\rm (Theorem 3 in \cite{MyB-2017})} Let $\mathcal{C}$ be a rational algebraic plane curve defined by a pa\-ra\-me\-tri\-za\-tion $\mathcal{P}(t)$, with limit point $P_L$. Let $P_1,\ldots,P_n$ be the singular points of $\cal C$, with multiplicities $m_1,\ldots,m_n$ respectively. Let us assume that they are ordinary singularities and that $P_i\neq P_L$ for $i=1,\ldots,n$. Then, it holds that
$$T(s)=\prod_{i=1}^nH_{P_i}(s)^{m_i-1}.$$
\end{theorem}

\begin{corollary}\label{C-tfunct}
Let $\cal C$  be a rational plane curve such that all its singularities are ordinary. Let $\mathcal{P}(t)$ be a parametrization of $\cal C$ such that $P_L$ is regular. It holds that
$\deg(T)=(d-1)(d-2).$
\end{corollary}

\para

Theorem \ref{T-tfunct} and Corollary \ref{C-tfunct} hold under the assumption that the limit point, $P_L$, is regular. In Theorem \ref{T-tfunct3} we eliminate this condition and  generalize both results to the case that $P_L$ is a singularity. For this purpose, we state the following theorem  whose proof will be presented in Section \ref{S-proof}.

\para

\begin{theorem}\label{T-tfunct2} Let $\mathcal{C}$ be a rational algebraic plane curve defined by a proper pa\-ra\-me\-tri\-za\-tion $\mathcal{P}(t)$ with limit point $P_L$. Let $P_1,\ldots,P_n$ and $P_L$ be the singularities of $\mathcal{C}$ and suppose that all of them are ordinary. For each $P_i \,\,(i=1,\ldots,n)$, let $m_i$ be its multiplicity and $H_{P_i}$ its fibre function. In addition, let $m_L$ and $H_L$ be the multiplicity and the fibre function of $P_L$. Then,
$$T(s)=\prod_{i=1}^nH_{P_i}(s)^{m_i-1}H_L(s)^{m_L-1}.$$
\end{theorem}

\begin{corollary}\label{C-tfunct2} It holds that
$$\deg(T)=(d-1)(d-2)-m_H(m_L-1).$$
\end{corollary}

\noindent\textbf{Proof:} First we note that if $P_L$ is regular, then $m_L=1$ and the result  follows from Corollary \ref{C-tfunct}. If $P_L$ is singular, from Theorem \ref{T-tfunct2}, we deduce that
$$\deg(T)=\sum_{i=1}^n(m_i-1)\deg(H_{P_i})+(m_L-1)\deg(H_L),$$
where $P_1,\ldots,P_n,\,P_L$ are the singularities of $\mathcal{C}$, and $m_1,\ldots,m_n,m_L$ their respective multiplicities. From Corollary \ref{C-multiplicidad}, we get that $\deg(H_{P_i})=m_i$. Furthermore, in Section \ref{S-pto-limite}, we show  that $\deg(H_L)$ provides the visible multiplicity of $P_L$; that is, $\deg(H_L)=m_L-m_H$, where $m_H$ is the hidden multiplicity. Then,
$$\deg(T)=\sum_{i=1}^nm_i(m_i-1)+(m_L-m_H)(m_L-1).$$
Finally, the corollary follows using the genus formula (see Chapter 3 in \cite{SWP})  which, in this case, implies that
$$\sum_{i=1}^nm_i(m_i-1)+m_L(m_L-1)=(d-1)(d-2).$$
\hfill $\Box$

\para

  Note that Theorem \ref{T-tfunct2} also holds if the limit point is regular since, in this case, $m_L=1$ and the corresponding factor disappears. Hence, by combining Theorem \ref{T-tfunct} and Theorem \ref{T-tfunct2}, we can state the following result, which does not impose any condition concerning the limit point.

\para
\begin{theorem}\label{T-tfunct3} Let $\mathcal{C}$ be a rational algebraic plane curve defined by a para\-me\-tri\-zation $\mathcal{P}(t)$. Let $P_1,\ldots,P_n$ be the singular points of $\cal C$, with multiplicities $m_1,\ldots,m_n$ respectively. Let us assume that they are ordinary singularities. Then, it holds that
$$T(s)=\prod_{i=1}^nH_{P_i}(s)^{m_i-1}.$$
\end{theorem}
\para
\begin{remark} \label{R-interp-T}
Using Theorem \ref{T-tfunct3}, we get that:
\begin{enumerate}
\item If for each factor $H_{P_i}(s)^{m_i-1}$ it holds that $\deg(H_{P_i})=m_i,\,i=1,\ldots,n$, then the limit point is regular.
\item If there is a factor $H_{P_{i_0}}(s)^{m_{i_0}-1}$ such that $\deg(H_{P_{i_0}})<m_{i_0}$, then $P_{i_0}$ is the limit point and $m_H=m_{i_0}-\deg(H_{P_{i_0}})$.
\end{enumerate}
\end{remark}

\begin{remark} \label{R-conjug-points}
In general, different conjugate roots of the T--function appear all together under a unique irreducible polynomial. These roots are associated to families of conjugated parametric points (see Definition 4 in \cite{MyB-2017}). In \cite{JSC-Perez} (Theorem 16), it is  shown that all the points in such a family have the same multiplicity and character. Let us assume  that $T(s)$ includes a factor $m(s)^{k-1}$, where $m(s)$ is an irreducible polynomial of degree $l$. Then, $m(s)$ contains the fibre functions of $l/k$ singular points of multiplicity $k$ (see Theorem 5 in \cite{MyB-2017}).
\end{remark}

\para

The following example shows how useful the above results are in order to study the singularities of a rational plane curve defined parametrically.

\para

\begin{example} In Example \ref{Ex-PL-reachable}, we consider the rational plane curve  $\mathcal{C}$  defined over $\Bbb C$ by the projective parametrization $\mathcal{P}(t)=$  $$(t^6+2t^5+3t^4+3t^2+4t^3+2t+1: -t^4-t^3-t^2-2t-1: t^7+3t^5+t^4+3t^3+t+2t^2+1).$$
We compute the T--function by applying (\ref{Eq-tfunct}), and we get that up to constants in $\Bbb C$, $$T(s)=(s^2-s-1)(s^3+s+1)^2$$
$$(s^8+3s^7+13s^6+22s^5+43s^4+47s^3+47s^2+44s+16)(s^2+1)^6(s+1).$$
In the following, we analyze each of the five factors of $T(s)$:

\begin{itemize}
\item The first factor is $f_1(s)=s^2-s-1$. It has degree 2 and its power is 1, which means that $f_1$ provides a double point. Indeed, $f_1$ has two complex conjugated roots, $1/2+1/2\sqrt{5}$ and $1/2-1/2\sqrt{5}$. By substituting them into the parametrization we get the affine double point $P_1=(1:-1/5:1)$.
\item The second factor is $f_2(s)=(s^3+s+1)^2$. It has degree 3 and  its power is 2, which means that $f_2$ provides a triple point. If we compute the three roots of $f_2$ and substitute them in the parametrization we get the infinity triple point $P_2=(1:0:0)$.
\item The third factor is $f_3(s)=s^8+3s^7+13s^6+22s^5+43s^4+47s^3+47s^2+44s+16$. It is an irreducible polynomial of degree 8 and its power is 1. This implies that $f_3$ is associated to a set of double points which defines a family of conjugated parametric points (see Remark \ref{R-conjug-points}). More precisely, since the degree of the irreducible polynomial is 8, we deduce that this family contains 4 double points.
\item The fourth factor is $f_4(s)=(s^2+1)^6$. This factor may be difficult to interpret since one could think that the fibre function is $(s^2+1)$ and its power is 6. However, this is not correct. Actually, the fibre function is $(s^2+1)^2$ and its power is 3, which means that $f_4$ provides a   point of multiplicity $4$. In order to avoid mistakes, one has to compute one of the roots, namely $s=I$, the corresponding point of the curve, $P_4=\mathcal{P}(I)=(0:1:0)$, and the associate fibre function. From Remark \ref{R-afin-noafin}, we get that $H_{P_4}(s)=(s^2+1)^2$.
\item The last factor is $f_5(s)=s+1$, and it provides the point $P_5=\mathcal{P}(-1)=(0:0:1)$. The fibre function of $P_5$ is $H_{P_5}(s)=s+1$. Observe that $H_{P_5}$ has degree 1 and its power is 1. This situation is described in statement 2 of Remark \ref{R-interp-T}. More precisely,  we have that $m_5=2>\deg(H_{P_5})=1$, which implies that $P_5$ is the limit point. Its total multiplicity is $m_L=2$ and its hidden multiplicity is $m_H=1$.
\end{itemize}
In Figure \ref{F-pto-limite}, we plot the curve $\mathcal{C}$. One may see the two real affine double points $P_1$ and $P_5$.

\begin{figure}[h]
$$
\psfig{figure=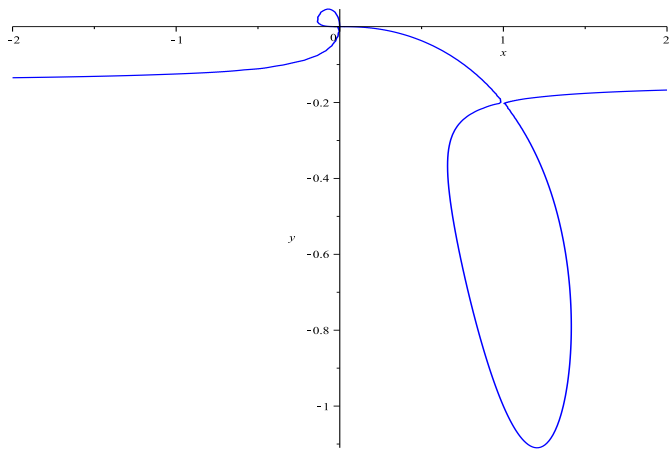,width=9cm,height=6cm,angle=0}
$$ \caption{The curve $\mathcal{C}$ has two real affine double points at $(1:-1/5:1)$ and $(0:0:1)$}\label{F-pto-limite}
\end{figure}

\end{example}

\subsection{The case of  rational space curves}\label{Sub-SpaceCurves}

Up to now,  we have dealt with rational plane curves, defined by para\-me\-tri\-za\-tions of the form $(p_1(t):p_2(t):p(t))\in\mathbb{P}^2(\mathbb{K}(t))$. Now, let $\mathcal{C}$ be a rational space curve defined by the projective proper parametrization $$\mathcal{P}(t)=(p_1(t):\cdots: p_n(t):p(t))\in\mathbb{P}^n(\mathbb{K}(t)),$$ where $\gcd(p_1,\ldots,p_n,p)=1$. In Section 4 in \cite{MyB-2017}, under the assumption that the limit point is regular,   we construct a polynomial, $T_E(s)$, which is equivalent to the T--function  introduced for plane curves,  and we prove that this polynomial  describes totally the singularities of $\cal C$,   since the factorization of  $T_E(s)$  provides the fibre functions of each singularity as well as its corresponding multiplicity. The idea is to construct a plane curve $\widehat{\mathcal{C}}$, defined over ${\Bbb K}(Z)$, where $Z=(Z_1,\ldots,Z_{n-2})$ and $Z_1,\ldots,Z_{n-2}$ are new variables, which contains  the information about the singularities of $\mathcal{C}$ and their multiplicities. Then, Theorem \ref{T-tfunct} may be used to get that information from $\widehat{\mathcal{C}}$.

\para

In the following we show that this idea may also be used for generalizing Theorem \ref{T-tfunct3} to the case of space curves. This will allow us to compute the singularities of any rational space  curve, from its parametric expression, even if one of those singularities is the limit point of the parametrization.

\para

\noindent For this purpose, let $\widehat{\mathcal{C}}$ be the plane curve defined by the parametrization
$$\widehat{\mathcal{P}}(t)=(\widehat{p}_1(t):\widehat{p}_2(t):\widehat{p}(t))=$$$$=(p_1(t):p_2(t)+Z_1p_{3}(t)+\cdots +Z_{n-2}p_n(t):p(t))\in\mathbb{P}^2((\mathbb{K}(Z))(t)).$$
We use this notation for the sake of simplicity but, we note that $\widehat{\mathcal{P}}(t)$ depends on $Z$. Observe that $\widehat{\mathcal{P}}(t)$ is a proper parametrization of $\widehat{{\cal C}}$ defined over the algebraic closure of  ${\Bbb K}(Z)$. In addition, let $\widehat{G}_1$, $\widehat{G}_2$ and $\widehat{G}_3$ be the equivalent polynomials to $G_1$, $G_2$ and $G_3$ (defined in (\ref{Eq-fibra-generica})), but constructed from the parametrization $\widehat{\mathcal{P}}(t)$. Similarly, let $\widehat{\delta}_i:=\degree_t(\widehat{G}_i)$, $\widehat{\lambda}_{ij}:=\min\{\widehat{\delta}_i,\widehat{\delta}_j\},\,i,j=1,2,3,\,i<j$, $$\widehat{G}_i^*(s,t):=\displaystyle\frac{\widehat{G}_i(s,t)}{t-s}\in ({\Bbb K}[Z])[s,t],\,\,\,i=1,2,3,$$
and
$$\widehat{R}_{ij}(s):=\Res_t(\widehat{G}_i^*,\widehat{G}_j^*)\in ({\Bbb K}[Z])[s],\,\,\, i,j=1,2,3,\,i<j.$$

\para

\noindent
The T--function of the parametrization  $\widehat{\mathcal{P}}(t)$
is given by
$$\widehat{T}(s)=\widehat{R}_{12}(s)/\widehat{p}(s)^{\widehat{\lambda}_{12}-1}.$$
It holds that $\widehat{T}(s)\in ({\Bbb K}[Z])[s]$  and
$$\widehat{T}(s)=\frac{\widehat{R}_{12}(s)}{\widehat{p}(s)^{\widehat{\lambda}_{12}-1}}=\frac{\widehat{R}_{13}(s)}{\widehat{p}_1(s)^{\widehat{\lambda}_{13}-1}}=\frac{\widehat{R}_{23}(s)}{\widehat{p}_2(s)^{\widehat{\lambda}_{23}-1}}$$
 (see \cite{MyB-2017}). Finally, we define the polynomial $$T_E(s)=\Content_Z\left(\widehat{T}(s)\right)\in {\Bbb K}[s].$$

The following theorem shows that $T_E(s)$ describes totally the singularities of $\cal C$,   since its factorization provides the fibre functions of the different singularities as well as their multiplicities. The result provided is similar to Theorem 6 in \cite{MyB-2017}, but the restriction for the limit point to be regular has been removed.

\para

Similarly as in  \cite{MyB-2017}, we need to impose the condition for the curve $\cal C$ not to have two or more {\it bad points} that is, points of the form $(0:a_2: a_3:\cdots:a_n: 0)$. Otherwise, we apply a change of coordinates, and we consider the new curve $\mathcal{C}^*$ parametrized by
$\mathcal{P}^*(t)=(p^*_1(t):p_2(t):\cdots: p_n(t):p(t))$, where $p^*_1=\sum_{i=1}^n\lambda_ip_i,\,\,\lambda_i\in {\Bbb K}$. By appropriately choosing the coefficients $\lambda_1,\ldots,\lambda_n$, we get that $\gcd(p^*_1,p)=1$ and  thus,  $\mathcal{C}^*$, does not have bad points. It is easy to check that, under these conditions, for each point $P=(a_1:a_2: a_3:\cdots:a_n: a_{n+1})\in {\cal C}$ we have another point $\widehat{P}=(a_1:a_2+Z_1a_{3}+\cdots +Z_{n-2}a_n: a_{n+1}) \in \widehat{{\cal C}}$ and this correspondence is bijective.

\para

\begin{theorem}\label{T-tfunct-spa}
Let $\mathcal{C}$ be a rational algebraic space curve, defined by a parametrization $\mathcal{P}(t)$. Let $P_1,\ldots,P_n$ be the singular points of $\cal C$, with multiplicities $m_1,\ldots,m_n$ respectively. Let us assume that they are ordinary singularities. Then, it holds that
$$T_E(s)=\prod_{i=1}^nH_{P_i}(s)^{m_i-1}.$$
\end{theorem}

The theorem may be proved similarly to Theorem 6 in \cite{MyB-2017}. That proof was based on the fact that $\widehat{P}=(a_1:a_2+Z_1a_{3}+\cdots +Z_{n-2}a_n: a_{n+1}), \,a_i\in {\Bbb K},\,i=1,\ldots,n+1,$ is a singularity of $\widehat{{\cal C}}$ of multiplicity $m$ if and only if $P=(a_1:a_2: a_3:\cdots:a_n: a_{n+1})$ is a singularity of ${\cal C}$ of multiplicity $m$, and that $H_{\widehat{P}}(s)=H_P(s)$. Now, we have just to observe that both statements also hold for the limit points $\widehat{P}_L$ and $P_L$.

\para

We observe that $\widehat{{\cal C}}$ may have additional singularities which fiber is in the algebraic closure of $\mathbb{K}(Z)\setminus\mathbb{K}$. These points do not have a correspondence with any point of $\cal C$. In order to remove them, we consider the content w.r.t. $Z$ of the polynomial $\widehat{T}(s)$ (we recall that $T_E(s)=\Content_Z\left(\widehat{T}(s)\right)\in {\Bbb K}[s]$).

\para

Note that the polynomial $H_{P}$  represents the fibre function of a point $P$ in the space curve $\cal C$ computed from $\mathcal{P}(t)$; i.e. the roots of $H_{P}$ are the fibre of    $P\in {\cal C}$ (this notion was introduced in Definition \ref{D-fibre-func} for a given plane curve but it can be easily generalized for space curves).

\para

We recall that from the fibre function of a point, one may determine the multiplicity of the point  as well as its fibre and the tangent lines at that point (see Section \ref{S-notacion}). The method presented generalizes the results obtained in \cite{Rubio}, since a complete classification of the singularities of a given space curve, via the factorization of a resultant, is obtained.

\para

\begin{example}
Let $\mathcal{C}$ be the rational space curve defined by the projective parametrization
$\mathcal{P}(t)=(p_1(t):p_2(t):p_3(t):p(t))\in\mathbb{P}^3(\mathbb{C}(t)),$
where
$$\begin{array}{l}
    p_1(t)=t^5+t^4-16t^3+16t^2-17t+15 \\
    p_2(t)=-t^5+4t^4-4t^3+3t^2+t-3 \\
    p_3(t)=t^5+t^4-19t^3+13t^2+34t-30\\
    p(t)=t^6-t^5+12t^4+23t^3+40t^2+27t+15.
  \end{array}$$
We consider the plane curve $$\widehat{\mathcal{P}}(t)=(p_1(t):p_2(t)+Zp_{3}(t):p(t))\in\mathbb{P}^2((\mathbb{C}(Z))(t))$$
and we compute the corresponding T--function which is given, up to constants in $\Bbb C$, as
$$\widehat{T}(s)=(t-1)^2(t-3)^2L(s,Z),$$
where $L(s,Z)$ is an irreducible polynomial which depends on $Z$. Observe that $L(s,Z)$ provides the singularities of the plane curve which fiber is in the algebraic closure of $\mathbb{K}(Z)\setminus\mathbb{K}$. As we stated above, these points have to be removed, so we eliminate this factor by computing the content w.r.t. $Z$. Thus, we get
$$T_E(s)=(t-1)^2(t-3)^2.$$
Since the power of the factors $(t-1)$ and $(t-3)$ is 2, we deduce that they are associated to the triple point  $P=\mathcal{P}(1)=\mathcal{P}(3)=(0:0:0:1)$. However, the fibre function, $H_P(t)=(t-1)(t-3)$, has degree 2, which does not agree with the multiplicity. Hence, we deduce that $P$ is the limit point of the parametrization and part of its multiplicity is not reflected by its fibre (see Section \ref{S-pto-limite}).

\para

In Figure \ref{F-espacial}, we plot the curve  $\mathcal{C}$. Observe that $P$ is reached by ${\mathcal P}(t)$ just two times (note that we are plotting $P$ from $\cP(t)$ with $-200\leq t\leq 200$). In the limit, when $t$ tends to $\pm\infty$, this point  would be reached a third time.
\begin{figure}[h]
$$
\psfig{figure=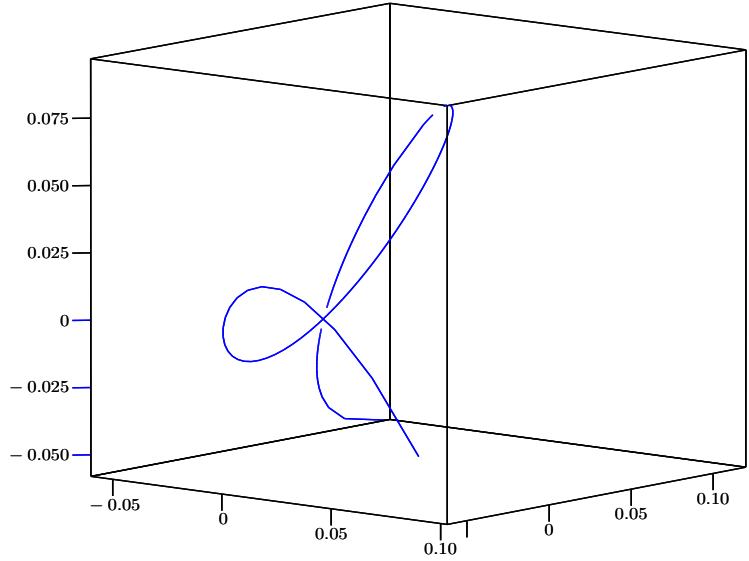,width=12cm,height=10cm}
$$\caption{Curve $\mathcal{C}$ plotted in a neighborhood of the triple point $P$}\label{F-espacial}
\end{figure}
\end{example}

\section{Proof  of  Theorem \ref{T-tfunct2} in Section \ref{S-tfunct}}\label{S-proof}

This section is devoted to prove  Theorem \ref{T-tfunct2} stated in Section \ref{S-tfunct}. For this purpose, throughout this section, we  use the reparametrization
\begin{equation}\label{Eq-Q}\mathcal{Q}(t):=\mathcal{P}\left(\frac{\theta t}{t-1}\right)\end{equation}
introduced in Remark \ref{R-raiz-cero}.  Observe that
$$\mathcal{Q}(t)=\left(p_1\left(\frac{\theta t}{t-1}\right)(t-1)^d:p_2\left(\frac{\theta t}{t-1}\right)(t-1)^d:p\left(\frac{\theta t}{t-1}\right)(t-1)^d\right)$$
and then, we may write $\mathcal{Q}(t)=(q_1(t): q_2(t): q(t))$, where
\begin{equation}\label{Eq-Q-componentes}\left\{\begin{array}{r}q_1(t)=a_0(t-1)^d+a_1\theta t(t-1)^{d-1}+\cdots+a_{d}\theta^dt^d\\q_2(t)=b_0(t-1)^d+b_1\theta t(t-1)^{d-1}+\cdots+b_{d}\theta^dt^d\\
q(t)=c_0(t-1)^d+c_1\theta t(t-1)^{d-1}+\cdots+c_{d}\theta^dt^d.\end{array}\right.\end{equation}
Note that, with this parametrization, $P_L$ is not the limit point but a conventional point which can be obtained as $P_L=\mathcal{Q}(1)$. The new limit point is
\begin{equation}\label{Eq-QL}Q_L=\lim_{t\rightarrow\infty}\mathcal{Q}(t)/t^d=\mathcal{P}(\theta),\end{equation}
and  we can choose $\theta\in\mathbb{K}\setminus\{0\}$ so that $Q_L$ is not a singularity. In addition, if $p(\theta)\neq 0$, we get that $Q_L$ is an affine regular point. Note that every point of the curve reached via $\mathcal{P}(t)$, but $Q_L$, is also reached via $\mathcal{Q}(t)$. Indeed, for any $s_0\neq\theta$ it holds that $\mathcal{Q}\left(\frac{s_0}{s_0-\theta}\right)=\mathcal{P}(s_0)$. The only exception arises when $s_0=\theta$, since $\mathcal{P}(\theta)=Q_L$.

\para

Taking  into account this last statement, we have that for each $s_i\in\mathcal{F}_{\mathcal{P}}(P_L)$ there exists  $t_i=s_i/(s_i-\theta)\in\mathcal{F}_{\mathcal{Q}}(P_L)$ (note that $s_i\not=\theta$ since $\theta\notin\mathcal{F}_{\mathcal{P}}(P_L)$; we recall that $\mathcal{P}(\theta)=Q_L\not=P_L$). However, there is a value in $\mathcal{F}_{\mathcal{Q}}(P_L)$ which does not have a correspondence in $\mathcal{F}_{\mathcal{P}}(P_L)$; namely $t=1$. Note that this value belongs to $\mathcal{F}_{\mathcal{Q}}(P_L)$, since $\mathcal{Q}(1)=P_L$, but there is no $s_i\in \mathbb{K}$ such that $\displaystyle\frac{s_i}{s_i-\theta}=1$ (note that $\theta\in\mathbb{K}\setminus\{0\}$). As a consequence, the fibre function of $P_L$ under the parametrization $\mathcal{Q}(t)$ is
\begin{equation}\label{Eq-HLQ2}\widetilde{H}_L(t)=(t-1)^r\prod_{i=1}^n\left(t-\frac{s_i}{s_i-\theta}\right)^{k_i},\end{equation}
where $k_1+\cdots+k_n$ is the visible multiplicity and $r$ is the hidden multiplicity (see the analogy with (\ref{Eq-HLQ})). Hence, we observe that $\deg(\widetilde{H}_L)=k_1+\cdots+k_n+r=\deg(H_L)+m_H=m_L$.

\para

Summarizing, we have a new parametrization of $\cal C$ such that the limit point, $Q_L$, is regular; this allows us to apply the results obtained in \cite{MyB-2017}. For this purpose, let  $\widetilde{G}_1(s,t)$, $\widetilde{G}_2(s,t)$ and $\widetilde{G}_3(s,t)$ be the equivalent polynomials to $G_1(s,t)$, $G_2(s,t)$ and $G_3(s,t)$ (see (\ref{Eq-fibra-generica})) computed from the new parametrization ${\cal Q}(t)$. That is,
\begin{equation}\label{Eq-Gtilde}\left\{\begin{array}{l}\widetilde{G}_1(s,t)=q_1(s)q(t)-q(s)q_1(t)\\\widetilde{G}_2(s,t)=q_2(s)q(t)-q(s)q_2(t)\\\widetilde{G}_3(s,t)=q_1(s)q_2(t)-q_2(s)q_1(t).\end{array}\right.\end{equation}
In addition, let $\widetilde{\delta}_i=\degree_t(\widetilde{G}_i)$, $\widetilde{\lambda}_{ij}=\min\{\widetilde{\delta}_i,\widetilde{\delta}_j\},\,i,j=1,2,3$, and
\begin{equation}\label{Eq-Rtilde}\widetilde{R}_{ij}(s)=\Res_t\left(\frac{\widetilde{G}_i(s,t)}{t-s},\frac{\widetilde{G}_j(s,t)}{t-s}\right),\quad \mbox{for $i,j=1,2,3$,\,$i < j$.}\end{equation}

\para

\noindent The T--function obtained from $\mathcal{Q}$ is
\begin{equation}\label{Eq-DefTQ}\widetilde{T}=\widetilde{R}_{12}/q^{\widetilde{\lambda}_{12}-1}=\widetilde{R}_{13}/q_1^{\widetilde{\lambda}_{13}-1}=\widetilde{R}_{23}/q_2^{\widetilde{\lambda}_{23}-1}.\end{equation}

\noindent Now, Theorem \ref{T-tfunct} and Corollary \ref{C-tfunct} state respectively that
$$\widetilde{T}(s)=\prod_{i=1}^n\widetilde{H}_{P_i}(s)^{m_i-1}\widetilde{H}_L(s)^{m_L-1}$$
and that
$$\deg(\widetilde{T})=(d-1)(d-2).$$

\para

We will use both statements later but, before, let us introduce the following technical lemma, which describes the relation between $R_{12}$ and $\widetilde{R}_{12}$ under the assumption that $\deg_t(G_1)=\deg_t(G_2)$.

\begin{lemma}\label{L-RP-RQ}
Let $\deg_t(G_1)=\deg_t(G_2)$. Then, it holds that
$$\widetilde{R}_{12}(s)=(s-1)^{2(d-1)^2}R_{12}\left(\frac{\theta s}{s-1}\right).$$
\end{lemma}

\noindent\textbf{Proof:} We prove the lemma by considering two steps. In the first one, we obtain the form of $R_{12}\left(\frac{\theta s}{s-1}\right)$. In the second one, we compute $\widetilde{R}_{12}(s)$ and we compare   it with  $R_{12}\left(\frac{\theta s}{s-1}\right)$.

\begin{center}Step 1\end{center}

\noindent
First, we recall that
$$R_{12}(s)=\Res_t\left(\frac{G_1(s,t)}{t-s},\frac{G_2(s,t)}{t-s}\right),$$
where $G_1(s,t)=p_1(s)p(t)-p(s)p_1(t)$ and
$G_2(s,t)=p_2(s)p(t)-p(s)p_2(t)$. In addition, since $\delta_1=\delta_2$, where  $\delta_1:=\deg_t(G_1)$ and $\delta_2:=\deg_t(G_2)$, and taking into account that $\delta_1=\max\{d_1,d_3\}$ and $\delta_2=\max\{d_2,d_3\}$ (see Remark 1 in \cite{MyB-2017}), we get that $\delta_1=\delta_2=\max\{d_1,d_2,d_3\}=d$.

\para

In the following, we denote $G_1^*:=G_1/(t-s)$ and $G_2^*:=G_2/(t-s)$, and thus
$R_{12}(s)=\Res_t(G^*_1,G^*_2)$, where $\deg_t(G_1^*)=\deg_t(G_2^*)=d-1$ (note that $\deg_s(G_1^*)=\deg_s(G_2^*)=d-1$). If we see $G_1^*$ and $G_2^*$ as polynomials in the variable $t$, that is,
$G_1^*,G_2^*\in(\mathbb{K}[s])[t]$, we may write
\begin{equation}\label{Eq-factor-G}G^*_1(s,t)=\lc_t(G^*_1)\prod_{i=1}^{d-1}(t-\alpha_i(s))\text{ and } G^*_2(s,t)=\lc_t(G^*_2)\prod_{j=1}^{d-1}(t-\beta_j(s)),\end{equation}
where $\lc_t(G^*_1)$ and $\lc_t(G^*_2)$ are their respective leader coefficients in ${\Bbb K}[s]$ and
$\alpha_1(s),\ldots,\alpha_{d-1}(s)$ and
$\beta_1(s),\ldots,\beta_{d-1}(s)$ their $d-1$
roots, respectively (that is,
$G_1^*(s,\alpha_i(s))=G_2^*(s,\beta_j(s))=0$ for
$i,j=1,\ldots,d-1$). Now, taking into account the properties of the resultant, we have that
\begin{equation}\label{Eq-resultante-P}R_{12}(s)=\lc_t(G^*_1)^{d-1}\lc_t(G^*_2)^{d-1}\prod_{i=1}^{d-1}\prod_{j=1}^{d-1}(\alpha_i(s)-\beta_j(s)).\end{equation}
By substituting $t=\theta$ in (\ref{Eq-factor-G}), we obtain
\begin{equation}\label{Eq-factor-G-2}G^*_1(s,\theta)=\lc_t(G^*_1)\prod_{i=1}^{d-1}(\theta-\alpha_i(s))\text{ and } G^*_2(s,\theta)=\lc_t(G^*_2)\prod_{j=1}^{d-1}(\theta-\beta_j(s))\end{equation}
so, up to constants in $\mathbb{K}\setminus\{0\}$,
\begin{equation}\label{Eq-lc-GP}\lc_t(G^*_1)=\frac{G^*_1(s,\theta)}{\prod_{i=1}^{d-1}(\alpha_i(s)-\theta)} \quad\text{ and }\quad \lc_t(G^*_2)=\frac{G^*_2(s,\theta)}{\prod_{j=1}^{d-1}(\beta_j(s)-\theta)}.\end{equation}

In the following, we denote $\Phi_1(s):=\prod_{i=1}^{d-1}(\alpha_i(s)-\theta)$ and $\Phi_2(s):=\prod_{j=1}^{d-1}(\beta_j(s)-\theta)$. Thus, taking into account that $G_i^*(s,t):=G_i(s,t)/(t-s)$, we get that
\begin{equation}\label{Eq-lc-GP-2}\lc_t(G^*_1)=\frac{G_1(s,\theta)}{(s-\theta)\Phi_1(s)} \quad\text{ and }\quad \lc_t(G^*_2)=\frac{G_2(s,\theta)}{(s-\theta)\Phi_2(s)}.\end{equation}
By substituting both expressions in (\ref{Eq-resultante-P}), we obtain that
\begin{equation}\label{Eq-resultante-P-2}R_{12}(s)=\left(\frac{G_1(s,\theta)}{(s-\theta)\Phi_1(s)}\right)^{d-1}\left(\frac{G_2(s,\theta)}{(s-\theta)\Phi_2(s)}\right)^{d-1}\prod_{i=1}^{d-1}\prod_{j=1}^{d-1}(\alpha_i(s)-\beta_j(s))\end{equation}
and therefore, $$\displaystyle R_{12}\left(\frac{\theta s}{s-1}\right)=\left(\frac{G_1\left(\frac{\theta s}{s-1},\theta\right)}{\left(\frac{\theta s}{s-1}-\theta\right)\Phi_1\left(\frac{\theta s}{s-1}\right)}\right)^{d-1}$$$$\left(\frac{G_2\left(\frac{\theta s}{s-1},\theta\right)}{\left(\frac{\theta s}{s-1}-\theta\right)\Phi_2\left(\frac{\theta s}{s-1}\right)}\right)^{d-1}\prod_{i=1}^{d-1}\prod_{j=1}^{d-1}\left(\alpha_i\left(\frac{\theta s}{s-1}\right)-\beta_j\left(\frac{\theta s}{s-1}\right)\right)$$
\begin{equation}\label{Eq-resultante-P-3}=\left(\frac{G_1\left(\frac{\theta s}{s-1},\theta\right)}{\Phi_1\left(\frac{\theta s}{s-1}\right)}\right)^{d-1}\left(\frac{G_2\left(\frac{\theta s}{s-1},\theta\right)}{\Phi_2\left(\frac{\theta s}{s-1}\right)}\right)^{d-1}\left(\frac{s-1}{\theta}\right)^{2(d-1)}$$$$\prod_{i=1}^{d-1}\prod_{j=1}^{d-1}\left(\alpha_i\left(\frac{\theta s}{s-1}\right)-\beta_j\left(\frac{\theta s}{s-1}\right)\right).\end{equation}

\begin{center}Step 2\end{center}

First, we observe that the coefficient of the term $t^d$ in $q(t)$ is $p(\theta)$. Thus,  $\deg(q)=d$ (note that $\theta$ is such that $p(\theta)\neq 0$) and hence,
$\deg_t(\widetilde{G}_1)=d$. Indeed:  we have that $\widetilde{G}_1(s,t)=q_1(s)q(t)-q(s)q_1(t)$ (see (\ref{Eq-Gtilde})), so it follows   that $\deg_t(\widetilde{G}_1)=d$ if $\deg(q_1)<d$. On the other hand, if $\deg(q_1)=d$, it could happen that $\deg_t(\widetilde{G}_1)<d$  if $q_1(s)q(\theta)-q(s)q_1(\theta)=0$ but this would imply that $q_1(s)/q(s)$ is a constant and, thus, $\mathcal{C}$ is a line, which is impossible by assumption. Reasoning similarly, we deduce that $\deg_t(\widetilde{G}_2)=d$.

\para

Let $\widetilde{G}_1^*:=\widetilde{G}_1/(t-s)$ and
$\widetilde{G}_2^*:=\widetilde{G}_2/(t-s)$. Then,
$\widetilde{R}_{12}(s)=\Res_t(\widetilde{G}^*_1,\widetilde{G}^*_2)$. Similarly as before, we have that, $\widetilde{G}_1^*$ and $\widetilde{G}_2^*$ are polynomials
in the variables $s$ and $t$, with degree $d-1$ in both variables. Moreover, reasoning as in Step 1, we have that
\begin{equation}\label{Eq-factor-GQ}\widetilde{G}^*_1(s,t)=\lc_t(\widetilde{G}^*_1)\prod_{i=1}^{d-1}(t-\widetilde{\alpha}_i(s)),\quad\quad
\widetilde{G}^*_2(s,t)=\lc_t(\widetilde{G}^*_2)\prod_{j=1}^{d-1}(t-\widetilde{\beta}_j(s)),\end{equation}
and
\begin{equation}\label{Eq-resultante-Q}\widetilde{R}_{12}(s)=\lc_t(\widetilde{G}^*_1)^{d-1}\lc_t(\widetilde{G}^*_2)^{d-1}\prod_{i=1}^{d-1}\prod_{j=1}^{d-1}(\widetilde{\alpha}_i(s)-\widetilde{\beta}_j(s)).\end{equation}
Note also that $\lc_t(\widetilde{G}^*_1)=\lc_t(\widetilde{G}_1)=$
$$q_1(s)(c_0+c_1\theta+\cdots +c_d\theta^d)-q(s)(a_0+a_1\theta+\cdots +a_d\theta^d)=q_1(s)p(\theta)-q(s)p_1(\theta)$$ which can be written as
$$p_1(\theta s/(s-1))(s-1)^dp(\theta)-p(\theta s/(s-1))(s-1)^dp_1(\theta)=G_1(\theta s/(s-1),\theta)(s-1)^d.$$
Reasoning similarly with  $\lc_t(\widetilde{G}^*_2)$,  we conclude that
\begin{equation}\label{Eq-lc-GQ}\lc_t(\widetilde{G}^*_1)=G_1\left(\frac{\theta s}{s-1},\theta\right)(s-1)^d \text{ and }
\lc_t(\widetilde{G}^*_2)=G_2\left(\frac{\theta s}{s-1},\theta\right)(s-1)^d.\end{equation}
Now, we focus on the roots
$\widetilde{\alpha}_1(s),\ldots,\widetilde{\alpha}_{d-1}(s)$. On the one side, all of them verify that
$\widetilde{G}_1^*(s,\widetilde{\alpha}_i(s))=0$. On the other side, note that
$$\widetilde{G}_1(s,t)=q_1(s)q(t)-q(s)q_1(t)=$$
$$=p_1\left(\frac{\theta s}{s-1}\right)(s-1)^dp\left(\frac{\theta t}{t-1}\right)(t-1)^d-p\left(\frac{\theta s}{s-1}\right)(s-1)^dp_1\left(\frac{\theta t}{t-1}\right)(t-1)^d$$
$$=(s-1)^d(t-1)^dG_1\left(\frac{\theta s}{s-1},\frac{\theta t}{t-1}\right).$$
Thus, $\widetilde{G}_1^*(s,\widetilde{\alpha}_i(s))=0$, which is equivalent to
$$G_1\left(\frac{\theta s}{s-1},\frac{\theta \widetilde{\alpha}_i(s)}{\widetilde{\alpha}_i(s)-1}\right)=0,$$ implies that
$$\alpha_i\left(\frac{\theta s}{s-1}\right)=\frac{\theta \widetilde{\alpha}_i(s)}{\widetilde{\alpha}_i(s)-1}$$ and thus $$\widetilde{\alpha}_i(s)=\frac{\alpha_i\left(\frac{\theta s}{s-1}\right)}{\alpha_i\left(\frac{\theta s}{s-1}\right)-\theta}.$$
Reasoning similarly with $\widetilde{\beta}_1(s),\ldots,\widetilde{\beta}_{d-1}(s)$,
we deduce that
\begin{equation}\label{Eq-aux-raices}\widetilde{\alpha}_i(s)=\frac{\alpha_i\left(\frac{\theta s}{s-1}\right)}{\alpha_i\left(\frac{\theta s}{s-1}\right)-\theta} \quad\text{ and
}\quad \widetilde{\beta}_j(s)=\frac{\beta_j\left(\frac{\theta s}{s-1}\right)}{\beta_j\left(\frac{\theta s}{s-1}\right)-\theta}\end{equation}
for each $i,j=1,\ldots,d-1$.

\para

\noindent Now, by substituting (\ref{Eq-lc-GQ}) and (\ref{Eq-aux-raices}) on
(\ref{Eq-resultante-Q}), we get
$$\widetilde{R}_{12}(s)=\left(G_1\left(\frac{\theta s}{s-1},\theta\right)(s-1)^d\right)^{d-1}\left(G_2\left(\frac{\theta s}{s-1},\theta\right)(s-1)^d\right)^{d-1}$$$$
\prod_{i=1}^{d-1}\prod_{j=1}^{d-1}\left(\frac{\alpha_i\left(\frac{\theta s}{s-1}\right)}{\alpha_i\left(\frac{\theta s}{s-1}\right)-\theta}-\frac{\beta_j\left(\frac{\theta s}{s-1}\right)}{\beta_j\left(\frac{\theta s}{s-1}\right)-\theta}\right)=$$
$$=G_1\left(\frac{\theta s}{s-1},\theta\right)^{d-1}G_2\left(\frac{\theta s}{s-1},\theta\right)^{d-1}(s-1)^{2d(d-1)}$$$$\prod_{i=1}^{d-1}\prod_{j=1}^{d-1}\left(\frac{\alpha_i\left(\frac{\theta s}{s-1}\right)\left(\beta_j\left(\frac{\theta s}{s-1}\right)-\theta\right)-\beta_j\left(\frac{\theta s}{s-1}\right)\left(\alpha_i\left(\frac{\theta s}{s-1}\right)-\theta\right)}{\left(\alpha_i\left(\frac{\theta s}{s-1}\right)-\theta\right)\left(\beta_j\left(\frac{\theta s}{s-1}\right)-\theta\right)}\right)=$$
$$=G_1\left(\frac{\theta s}{s-1},\theta\right)^{d-1}G_2\left(\frac{\theta s}{s-1},\theta\right)^{d-1}(s-1)^{2d(d-1)}$$$$\frac{\prod_{i=1}^{d-1}\prod_{j=1}^{d-1}\left(\alpha_i\left(\frac{\theta s}{s-1}\right)-\beta_j\left(\frac{\theta s}{s-1}\right)\right)\theta}{\prod_{i=1}^{d-1}\left(\alpha_i\left(\frac{\theta s}{s-1}\right)-\theta\right)^{d-1}\prod_{j=1}^{d-1}\left(\beta_j\left(\frac{\theta s}{s-1}\right)-\theta\right)^{d-1}}=$$

\para

$$=\left(\frac{G_1\left(\frac{\theta s}{s-1},\theta\right)}{\Phi_1\left(\frac{\theta s}{s-1}\right)}\right)^{d-1}\left(\frac{G_2\left(\frac{\theta s}{s-1},\theta\right)}{\Phi_2\left(\frac{\theta s}{s-1}\right)}\right)^{d-1}
(s-1)^{2d(d-1)}\theta^{(d-1)^2}$$$$\prod_{i=1}^{d-1}\prod_{j=1}^{d-1}\left(\alpha_i\left(\frac{\theta s}{s-1}\right)-\beta_j\left(\frac{\theta s}{s-1}\right)\right).$$

\noindent Finally, by comparing this expression with (\ref{Eq-resultante-P-3}),
we observe that, up to constants,
$$\frac{\widetilde{R}_{12}(s)}{(s-1)^{2d(d-1)}}=\frac{R_{12}(\theta s/(s-1))}{(s-1)^{2(d-1)}},$$
which proves the lemma.\hfill $\Box$

\para

\begin{remark}\label{R-RP-RQ}
Reasoning as in the proof of Lemma \ref{L-RP-RQ}, one may prove that
\begin{enumerate}
\item If $p_1(\theta)\neq 0$, then $\widetilde{R}_{13}(s)=(s-1)^{2(d-1)^2}R_{13}\left(\frac{\theta s}{s-1}\right).$
\item If $p_2(\theta)\neq 0$, then $\widetilde{R}_{23}(s)=(s-1)^{2(d-1)^2}R_{23}\left(\frac{\theta s}{s-1}\right).$
\end{enumerate}\end{remark}
\para

 Now,  we are ready to deal with the main item of the section. More precisely,  we prove Theorem \ref{T-tfunct2} (see Section \ref{S-tfunct}), where   Theorem \ref{T-tfunct} is generalized to the case that the limit point, $P_L$, is a singularity of the curve. In order to show this result, Lemma \ref{L-RP-RQ} and Remark \ref{R-RP-RQ} will be required.

\para
\para

\noindent\textbf{Proof of Theorem \ref{T-tfunct2}}

\para

\noindent We consider two steps in the proof of the theorem.  First, we assume that $\delta_1=\delta_2$, which allows us to use Lema \ref{L-RP-RQ} and  second,  we eliminate this requirement and we prove that the result holds anyway.

\begin{center}Step 1\end{center}

According to Lema \ref{L-tfunct},
for each singularity $P_i, i=1,\ldots,n$ it holds that
$T(s)=H_{P_i}(s)^{m_i-1}T_i^*(s)$, where $T_i^*$ is a polynomial such that $\gcd(H_{P_i},T_i^*)=1$. Note that
$\gcd(H_{P_i},H_{P_j})=1$ for $i\neq j$ (otherwise,
one single value of the parameter would generate two different points of the curve). Hence, we have that
\begin{equation}\label{Eq-Tfactor}T(s)=\prod_{i=1}^nH_{P_i}(s)^{m_i-1}V(s),\end{equation}
where $V(s)$ is a polynomial such that $\gcd(H_{P_i},V)=1$, for $i=1,\ldots,n$.

\para

Now, we focus on the polynomial $V$. We observe that if $V(s_0)=0$ then $T(s_0)=0$ which implies that $R_{12}(s_0)=R_{13}(s_0)=R_{23}(s_0)=0$ (see (\ref{Eq-tfunct})). Moreover, we note that $d=\max\{d_1,d_2,d_3\}$, so it can not happen that $a_d=b_d=c_d=0$. In the following we assume w.l.o.g., that $a_d\neq 0$ and we consider $R_{13}$ (if $a_d=0$ and $b_d\neq 0$, we would consider $R_{23}$ and if $a_d=b_d=0$ and $c_d\neq 0$ we would use $R_{12}$). Then, let $a_d\neq 0$. We have that if
$V(s_0)=0$, then $$R_{13}(s_0)=\Res_t(G^*_1(s,t),G^*_3(s,t))({s_0})=0,$$
and thus, one of the following statements hold:
\begin{enumerate}
\item The polynomials $G^*_1(s_0,t)$ and $G^*_3(s_0,t)$  have a common root, say $t=s_1\neq s_0$, which implies that the point $P=\mathcal{P}(s_0)=\mathcal{P}(s_1)$ is a singularity. However, it can not be $P=P_i$ for $i=1,\ldots,n$ since it would imply that $\gcd(H_{P_i},V)\neq 1$ (both polynomials would have $(t-s_0)$ as a common factor). Thus, we have that $\mathcal{P}(s_0)=P_L$ and  $H_L(s_0)=0$.
\item It holds that  $\gcd(\lc_t(G^*_1),\lc_t(G_3^*))({s_0})=0$. Note that
$$\lc_t(G^*_1)=\lc_t(G_1)=p_1(s)c_d-p(s)a_d=\phi_1^L(s) \mbox{ and} $$$$
\lc_t(G^*_3)=\lc_t(G_3)=p_1(s)b_d-p_2(s)a_d=\phi_3^L(s)$$ and then,
$\gcd(\lc_t(G^*_1),\lc_t(G_3^*))=\gcd(\phi_1^L(s),\phi_3^L(s))$. In addition, since we are assuming that $a_d\neq 0$, we may write
$$\phi_2^L(s)=\frac{b_d}{a_d}\phi_1^L(s)-\frac{c_d}{a_d}\phi_3^L(s),$$
and, hence,
$\gcd(\phi_1^L(s),\phi_3^L(s))=\gcd(\phi_1^L(s),\phi_2^L(s),\phi_3^L(s))=H_L(s)$. Finally we conclude that $$\gcd(\lc_t(G^*_1),\lc_t(G_3^*))=H_L(s).$$
\end{enumerate}
In both cases,   $V(s_0)=0$ implies that $H_L(s_0)=0$. Now, let us check that the reciprocal holds,  that is,   we assume that $H_L(s_0)=0$ (i.e., that $\mathcal{P}(s_0)=P_L$) and we   prove that $V(s_0)=0$.
By applying (\ref{Eq-resultante-P}) to $R_{13}$, we deduce that it is divided by $\gcd(\lc_t(G^*_1),\lc_t(G_3^*))=H_L(s)$. Thus, $H_L(s_0)=0$ implies that $R_{13}(s_0)=0$. In addition, by combining (\ref{Eq-tfunct}) and (\ref{Eq-Tfactor}),
$$R_{13}(s)=T(s)p_1(s)^{\lambda_{13}-1}=\prod_{i=1}^nH_{P_i}(s)^{m_i-1}V(s)p_1(s)^{\lambda_{13}-1}$$
so we get that
$$\prod_{i=1}^nH_{P_i}(s_0)^{m_i-1}V(s_0)p_1(s_0)^{\lambda_{13}-1}=0.$$
Observe that $\gcd(H_L,H_{P_i})=1$ (otherwise, we would have $P_i=P_L$). Hence, $\prod_{i=1}^nH_{P_i}(s_0)^{m_i-1}\neq 0$. On the other hand,
$\mathcal{P}(s_0)=P_L$ implies that $p_1(s_0)=a_d$ and since we are assuming $a_d\neq 0$,  we conclude that $V(s_0)=0$. Thus, we have proved that
$V(s_0)=0$ if and only if $H_L(s_0)=0$ and, as a consequence, we have that
$V(s)=H_L(s)^{\nu}$ for some $\nu\in\mathbb{N}$. By substituting above, we get
\begin{equation}\label{Eq-RQ}R_{13}(s)=\prod_{i=1}^nH_{P_i}(s)^{m_i-1}H_L(s)^{\nu}p_1(s)^{\lambda_{13}-1}.\end{equation}

\para

\noindent
Now,  we compute $\nu$. For this purpose, we use the reparametrization $$\mathcal{Q}(t)=\mathcal{P}\left(\frac{\theta t}{t-1}\right),$$ and we take $\theta$ such that $p_1(\theta)\neq 0$ and $\mathcal{P}(\theta)\neq P_i$, for $i=1,\ldots,n$. We consider the polynomials $\widetilde{G}_i,\,i=1,2,3$ and $\widetilde{R}_{ij}$ ($i,j=1,2,3$) introduced in (\ref{Eq-Gtilde}) and (\ref{Eq-Rtilde}), respectively. In addition, let $\widetilde{\delta}_i=\degree_t(\widetilde{G}_i)$ and $\widetilde{\lambda}_{ij}=\min\{\widetilde{\delta}_i,\widetilde{\delta}_j\}$ ($i,j=1,2,3$). Then, we may construct the T--function for the new parametrization, ${\cal Q}(t)$, as follows (see (\ref{Eq-DefTQ})):
$$\widetilde{T}=\widetilde{R}_{13}/q_1^{\widetilde{\lambda}_{13}-1}.$$

We recall that the limit point for ${\cal Q}(t)$ is $Q_L=\mathcal{P}(\theta)$ (see (\ref{Eq-QL})). Since we have chosen $\theta$ such that $\mathcal{P}(\theta)\neq P_i$ for $i=1,\ldots,n$, we are ensuring the new limit point to be regular. Thus, we can apply Corollary \ref{C-tfunct} and we get that $\deg(\widetilde{T})=(d-1)(d-2)$. In addition, we have that $p_1(\theta)\neq 0$ and thus $\deg(q_1)=d$, which implies that $\widetilde{\delta}_1=\widetilde{\delta}_3=\widetilde{\lambda}_{13}=d$. Hence,
$$\deg(\widetilde{R}_{13})=\deg(\widetilde{T}\cdot
q_1^{\widetilde{\lambda}_{13}-1})=(d-1)(d-2)+d(d-1)=2(d-1)^2.$$
Since $Q_L$ is a regular point, we may apply Theorem \ref{T-tfunct} and we get the following equality:
\begin{equation}\label{Eq-RQ-tilde}\widetilde{R}_{13}(s)=\prod_{i=1}^n\widetilde{H}_{P_i}(s)^{m_i-1}\widetilde{H}_L(s)^{m_L-1}q_1(s)^{d-1}\end{equation}
where $\widetilde{H}_{P_i}(s)$ is the fibre function of the singular point $P_i$ ($i=1,\ldots,n$) and $\widetilde{H}_{L}(s)$ is the fibre function of the singular point $P_L$ (note that now $P_L$ is not the limit point). All these fibre functions are computed from the parametrization $\mathcal{Q}(t)$.

\para

Now, let us compare the expressions in (\ref{Eq-RQ}) and (\ref{Eq-RQ-tilde}). By computing degrees on (\ref{Eq-RQ}), we have that
$$\deg_t(R_{13}(s))=\deg_t(\prod_{i=1}^nH_{P_i}(s)^{m_i-1})+\deg_t(H_L(s)^{\nu})+\deg_t(p_1(s)^{\lambda_{13}-1}).$$
On the other hand, from (\ref{Eq-RQ-tilde}),
$$\deg_t(\widetilde{R}_{13}(s))=\deg_t(\prod_{i=1}^n\widetilde{H}_{P_i}(s)^{m_i-1})+\deg_t(\widetilde{H}_L(s)^{m_L-1})+\deg_t(q_1(s)^{d-1}).$$
We observe that $P_i\neq P_L$ and $P_i\neq Q_L$ for each $i=1,\ldots,n$ (note that $Q_L=\mathcal{P}(\theta)$ and we have chosen $\theta$ such that $\mathcal{P}(\theta)\neq P_i$ for $i=1,\ldots,n$). Then, by applying Corollary \ref{C-multiplicidad}, we get that $\deg(H_{P_i}(s))=\deg(\widetilde{H}_{P_i}(s))=m_i$ for $i=1,\ldots,n$ and, therefore, $\deg(\prod_{i=1}^nH_{P_i}(s)^{m_i-1})=\deg(\prod_{i=1}^n\widetilde{H}_{P_i}(s)^{m_i-1})$. On the other hand, we have that $\deg(p_1)=d$ (we are assuming that $a_d\neq 0$) and $\deg(q_1)=d$ (the coefficient of $t^d$ in $q_1$ is $p(\theta)$ and we are assuming that $p(\theta)\neq 0$). Furthermore, $a_d\neq 0$ implies that $d_1=d\geq d_2,d_3$ and thus, $\delta_1=\delta_3=\lambda_{13}=d$, so we have that $\deg(p_1(s)^{\lambda_{13}-1})=\deg(q_1(s)^{d-1})$.

\para

\noindent
From the above statements, we deduce that
$$\deg(R_{13})-\deg(\widetilde{R}_{13})=\deg(H_L(s)^{\nu})-\deg(\widetilde{H}_L(s)^{m_L-1}).$$
In addition, we know from (\ref{Eq-HLQ2}) that $\deg(\widetilde{H}_L(s))$ provides the total multiplicity of $P_L$ (that is, $m_L$), while $\deg(H_L(s))$ provides just its visible multiplicity (that is, $m_L-m_H$). Thus, we have that
\begin{equation}\label{Eq-grado-RP1}\deg(R_{13})=2(d-1)^2+(m_L-m_H)\nu-m_L(m_L-1).\end{equation}

\para

In the following, we  compute $\deg(R_{13})$ in a different way and we compare the result obtained with (\ref{Eq-grado-RP1}). Since we are assuming $\delta_1=\delta_2$, we may use Lema \ref{L-RP-RQ} and Remark \ref{R-RP-RQ}; the last one states that
$$\widetilde{R}_{13}(s)=(s-1)^{2(d-1)^2}R_{13}\left(\frac{\theta s}{s-1}\right).$$
By applying the change $\displaystyle \frac{\theta s}{s-1}=t$ to the  above expression, we get
$$\widetilde{R}_{13}\left(\frac{t}{t-\theta}\right)=\left(\frac{t}{t-\theta}-1\right)^{2(d-1)^2}R_{13}(t)=\left(\frac{\theta}{t-\theta}\right)^{2(d-1)^2}R_{13}(t)$$
and hence,
\begin{equation}\label{Eq-R-Rtilde}R_{13}(t)=\left(\frac{t-\theta}{\theta}\right)^{2(d-1)^2}\widetilde{R}_{13}\left(\frac{t}{t-\theta}\right).\end{equation}

Now, let us analyze the factor $\widetilde{R}_{13}\left(\frac{t}{t-\theta}\right)$. According to (\ref{Eq-HLQ2}), we may write $\widetilde{H}_L(s)$ as
$\widetilde{H}_L(s)=(s-1)^{m_H}\overline{H}_L(s)$, where $m_H$ is the hidden multiplicity and $\overline{H}_L(s)$ is a polynomial such that $\overline{H}_L(1)\neq 0$. By substituting it on (\ref{Eq-RQ-tilde}),
we obtain:
$$\widetilde{R}_{13}(s)=\prod_{i=1}^n\widetilde{H}_{P_i}(s)^{m_i-1}\left((s-1)^{m_H}\overline{H}_L(s)\right)^{m_L-1}q_1(s)^{d-1}=$$
$$=(s-1)^{m_H(m_L-1)}\prod_{i=1}^n\widetilde{H}_{P_i}(s)^{m_i-1} \overline{H}_L(s)^{m_L-1}q_1(s)^{d-1}.$$
Note that $\overline{H}_L(1)\neq 0$ and $\widetilde{H}_{P_i}(1)\neq 0$ for any $i=1,\ldots,n$ ($\widetilde{H}_{P_i}(1)=0$ would imply that $P_i=\mathcal{Q}(1)=P_L$). In addition, from (\ref{Eq-Q-componentes}) we have that $q_1(1)=a_d\theta^d\neq 0$. Thus, we may write
$$\widetilde{R}_{13}(s)= (s-1)^{m_H(m_L-1)}\overline{R}_{13}(s),$$ where $\overline{R}_{13}(s)$ is a polynomial such that
$\overline{R}_{13}(1)\neq 0$. Hence,
$$\widetilde{R}_{13}\left(\frac{t}{t-\theta}\right)=\left(\frac{\theta}{t-\theta}\right)^{m_H(m_L-1)}\overline{R}_{13}\left(\frac{t}{t-\theta}\right)$$
and, by substituting in (\ref{Eq-R-Rtilde}), up to constants, we get that
$$R_{13}(t)=\left(\frac{t-\theta}{\theta}\right)^{2(d-1)^2}\left(\frac{\theta}{t-\theta}\right)^{m_H(m_L-1)}\overline{R}_{13}\left(\frac{t}{t-\theta}\right)=$$
$$=(t-\theta)^{2(d-1)^2-m_H(m_L-1)}\overline{R}_{13}\left(\frac{t}{t-\theta}\right).$$
Therefore, we deduce that
\begin{equation}\label{Eq-grado-RP2}\deg(R_{13})=2(d-1)^2-m_H(m_L-1)\end{equation}
and, by comparing (\ref{Eq-grado-RP1}) and (\ref{Eq-grado-RP2}),
we conclude that $\nu=m_L-1$.

\begin{center}Step 2\end{center}

In the following we prove that  condition $\delta_1=\delta_2$ is not necessary for the theorem to be hold. Let us assume w.l.o.g.
that $\delta_1:=\deg_t(G_1)>\delta_2:=\deg_t(G_2)$  that is,
$d=d_1>d_2,d_3$. Then, we consider the curve $\widehat{\mathcal{C}}$ defined by the parametrization
\begin{equation}\label{Eq-Pgorrito}\widehat{\mathcal{P}}(t)=(p_1(t):p_1(t)+\lambda p_2(t):p(t)),\end{equation}
where $\lambda\in {\Bbb K}$. Note that we are  applying a change of coordinates and then, for almost all values of $\lambda$, the curve  $\widehat{\cal C}$ has the same number of singularities that ${\cal C}$ and they are reached at the same values of the parameter $t\in {\Bbb K}$. That is, for each singularity $P_i$ in $\cal C$ there exists another singularity $\widehat{P}_i$ in $\widehat{\cal C}$ such that $H_{P_i}=H_{\widehat{P}_i}$, and reciprocally. Furthermore, since $P_i$ and $\widehat{P}_i$ have the same fibre function their multiplicities, say $m_i$ and $\widehat{m}_i$, are the same.

\para

\noindent The limit point for $\widehat{\mathcal{P}}(t)$ is
$$\widehat{P}_L=\lim_{t\rightarrow\infty}\frac{\widehat{\mathcal{P}}(t)}{t^d}=(a_d:a_d+\lambda b_d:c_d)$$
and it holds that $\widehat{\mathcal{P}}(t)=\widehat{P}_L$ if and only if $\mathcal{P}(t)=P_L$. Thus, $H_{\widehat{P}_L}(s)=H_{P_L}(s)$ and $\widehat{m}_L:=\mult_{\widehat{P}_L}(\widehat{\mathcal{C}})=m_L:=\mult_{P_L}(\mathcal{C})$.

\para

By appropriately choosing $\lambda$ in (\ref{Eq-Pgorrito}), we can force the singularities of the curve to keep their ordinary character. In addition, we can get that $a_d+\lambda b_d\neq 0$, which ensures that $\deg(p_1(t)+\lambda p_2(t))=\max\{d_1,d_2\}=d_1=d$. Now, for each $i=1,2,3$, let $\widehat{G}_i$ be the equivalent polynomial to $G_i$ but computed from $\widehat{\mathcal{P}}(t)$. Observe that $\widehat{G}_1=G_1$ and
$$\widehat{G}_2(s,t)=(p_1(s)+\lambda p_2(s))p(t)-p(s)(p_1(t)+\lambda p_2(t)),$$
which implies that $\deg_t(\widehat{G}_1)=\deg_t(\widehat{G}_2)=d$. Then, the conditions imposed in Step 1 are satisfied and the theorem holds for $\widehat{\mathcal{P}}(t)$. Thus,
$$\widehat{R}_{12}(s):=\Res_t\left(\frac{\widehat{G}_1}{t-s},\frac{\widehat{G}_2}{t-s}\right)=p(s)^{d-1}\left(\prod_{i=1}^n H_{\widehat{P}_i}(s)^{\widehat{m}_i-1}\right)H_{\widehat{P}_L}(s)^{\widehat{m}_L-1},$$
where $\widehat{P}_1,\ldots,\widehat{P}_n$ and $\widehat{P}_L$ are the singularities of $\widehat{\mathcal{C}}$, and $\widehat{m}_1,\ldots,\widehat{m}_n,\widehat{m}_L$ their respective multiplicities. On the other hand, we have remarked above that the singularities of ${\mathcal{C}}$ and $\widehat{\mathcal{C}}$ (and also their multiplicities) are the same. Then, we have that
\begin{equation}\label{Eq-Rgorrito1}\widehat{R}_{12}(s)=p(s)^{d-1}\prod_{i=1}^nH_{P_i}(s)^{m_i-1}H_L(s)^{m_L-1}.\end{equation}

\para

\noindent
Now we observe that  $$\widehat{G}_2(s,t)=(p_1(s)+\lambda p_2(s))p(t)-p(s)(p_1(t)+\lambda p_2(t))=$$$$p_1(s)p(t)-p(s)p_1(t)+\lambda (p_2(s)p(t)-p(s)p_2(t))=(G_1+\lambda G_2)(s,t).$$ Thus, by applying some well known properties of the resultants (see e.g. Appendix B in \cite{SWP}), we get that
$$\widehat{R}_{12}(s)=\Res_t\left(\frac{G_1}{t-s},\frac{G_1+G_2}{t-s}\right)=$$
$$=\lc_t\left(\frac{G_1}{t-s}\right)^{d-1}\prod_{i=1}^{d-1} \left(\frac{G_1(s,\alpha_i(s))+G_2(s,\alpha_i(s))}{\alpha_i(s)-s}\right)$$
where $\alpha_1(s),\ldots,\alpha_{d-1}(s)$ are the $d-1$ roots of the polynomial $G^*_1(s,t):=G_1(s,t)/(t-s)\in(\mathbb{K}[s])[t]$. That is, for each $i=1,\ldots,d-1$, we have that $\alpha_i(s)\neq s$ and $G_1(s,\alpha_i(s))=0$. Thus, the last formula may be written as
$$\widehat{R}_{12}(s)=\lc_t\left(\frac{G_1}{t-s}\right)^{d-\delta_2}\lc_t\left(\frac{G_1}{t-s}\right)^{\delta_2-1}\prod_{i=1}^{d-1} \left(\frac{G_2(s,\alpha_i(s))}{\alpha_i(s)-s}\right),$$
and, hence,
$$\widehat{R}_{12}(s)=\lc_t\left(\frac{G_1}{t-s}\right)^{d-\delta_2}\Res_t\left(\frac{G_1}{t-s},\frac{G_2}{t-s}\right).$$
Note that $\displaystyle \lc_t\left(\frac{G_1}{t-s}\right)=\lc_t(G_1)=p(s)$ since $d_3<d_1$, so we deduce that
\begin{equation}\label{Eq-Rgorrito2}\widehat{R}_{12}(s)=p(s)^{d-\delta_2}\Res_t\left(\frac{G_1}{t-s},\frac{G_2}{t-s}\right)=p(s)^{d-\delta_2}R_{12}(s).\end{equation}

\noindent Finally, by combining (\ref{Eq-Rgorrito1}) and (\ref{Eq-Rgorrito2}), we conclude that
$$R_{12}(s)=p(s)^{\delta_2-1}\prod_{i=1}^nH_{P_i}(s)^{m_i-1}H_L(s)^{m_L-1}.$$
The result follows from the assumption that $\delta_2=\min\{\delta_1,\delta_2\}=\lambda_{12}$.\hfill $\Box$

\para

\end{document}